\documentclass[11pt]{article}
\usepackage{amsfonts}
\usepackage{color}
\usepackage{amsmath,amssymb,comment}
\newtheorem{theo}{Theorem}[section]
\newtheorem{lem}[theo]{Lemma}

\newtheorem{defi}[theo]{Definition}

\newcommand{\mysection}[1]{\section{#1} \setcounter{equation}{0}}
\newcommand{\proof}{{\sc Proof.} \quad}
\newcommand{\proofc}{{\sc Proof} \ }
\newcommand{\be}{\begin{equation} \label}
\newcommand{\ee}{\end{equation}}
\newcommand{\bea}{\begin{eqnarray}\label}
\newcommand{\eea}{\end{eqnarray}}
\newcommand{\bas}{\begin{eqnarray*}}
\newcommand{\eas}{\end{eqnarray*}}
\newcommand{\bit}{\begin{itemize}}
\newcommand{\eit}{\end{itemize}}
\newcommand{\qed}{\hfill$\Box$ \vskip.2cm}
\newcommand{\nn}{\nonumber}
\newcommand{\R}{\mathbb{R}}
\newcommand{\N}{\mathbb{N}}
\newcommand{\pO}{\partial\Omega}

\newcommand{\eps}{\varepsilon}

\newcommand{\wto}{\rightharpoonup}
\newcommand{\wsto}{\stackrel{\star}{\rightharpoonup}}
\newcommand{\hra}{\hookrightarrow}
\newcommand{\io}{\int_\Omega}

\newcommand{\del}{\delta}
\newcommand{\al}{\alpha}
\newcommand{\lam}{\lambda}
\newcommand{\Lam}{\Lambda}

\newcommand{\bom}{\overline{\Omega}}
\newcommand{\Om}{\Omega}

\newcommand{\wh}{\widehat}

\newcommand{\hs}{\hspace*}

\newcommand{\vp}{\varphi}
\newcommand{\lbal}{\left\{ \begin{array}{l}}
\newcommand{\lball}{\left\{ \begin{array}{ll}}
\newcommand{\ear}{\end{array} \right.}

\newcommand{\abs}{\\[5pt]}

\newcommand{\adb}{\allowdisplaybreaks}
\newcommand{\tm}{T_{max}}

\newcommand{\ueps}{u_\eps}
\newcommand{\veps}{v_\eps}

\newcommand{\feps}{f_\eps}

\newcommand{\Teps}{\Theta_\eps}
\newcommand{\vepsx}{v_{\eps x}}
\newcommand{\vepsxx}{v_{\eps xx}}

\newcommand{\vepst}{v_{\eps t}}

\newcommand{\uepsxx}{u_{\eps xx}}

\newcommand{\uepst}{u_{\eps t}}
\newcommand{\Tepsx}{\Theta_{\eps x}}
\newcommand{\Tepsxx}{\Theta_{\eps xx}}
\newcommand{\Tepst}{\Theta_{\eps t}}

\begin{document}
\adb
%
%
\title{
Large-data regular solutions in a one-dimensional thermoviscoelastic
evolution problem involving temperature-dependent viscosities}
\author{
Michael Winkler\footnote{michael.winkler@math.uni-paderborn.de}\\
{\small Universit\"at Paderborn, Institut f\"ur Mathematik}\\
{\small 33098 Paderborn, Germany} }
\date{}
\maketitle
\begin{abstract}
\noindent 
The model
\bas
	\lbal
	u_{tt} = \big(\gamma(\Theta) u_{xt}\big)_x + au_{xx} - \big(f(\Theta)\big)_x, \\[1mm]
	\Theta_t = \Theta_{xx} + \gamma(\Theta) u_{xt}^2 - f(\Theta) u_{xt},
	\ear
\eas
for thermoviscoelastic evolution in one-dimensional Kelvin-Voigt materials is considered.
By means of an approach based on maximal Sobolev regularity theory of scalar parabolic equations,
it is shown that if $\gamma_0>0$ is fixed, then there exists $\del=\del(\gamma_0)>0$ with the property that
for suitably regular initial data of arbitrary size
an associated initial-boundary value problem posed in an open bounded interval admits a
global classical solution whenever $\gamma\in C^2([0,\infty))$ and $f\in C^2([0,\infty))$ are such that $f(0)=0$ and
$|f(\xi)| \le K_f \cdot (\xi+1)^\al$ for all $\xi\ge 0$ and some $K_f>0$ and $\al<\frac{3}{2}$, and that
\bas
	\gamma_0 \le \gamma(\xi) \le \gamma_0 + \del
	\qquad \mbox{for all } \xi\ge 0.
\eas
This is supplemented by a statement on global existence of certain strong solutions, particularly continuous
in both components, under weaker conditions on the initial data.\abs
\noindent {\bf Key words:} nonlinear acoustics; thermoviscoelasticity; viscous wave equation; maximal Sobolev regularity\\
{\bf MSC 2020:} 74H20 (primary); 74F05, 35L05, 35D30, 35B65 (secondary)
\end{abstract}
%
%
%
%
%
%
%
%
%
%
%
%
%
\newpage
\section{Introduction}\label{intro}
This manuscript is concerned with the problem
\be{00}
	\lbal
	u_{tt} = \big(\gamma(\Theta) u_{xt}\big)_x + au_{xx} - \big(f(\Theta)\big)_x, \\[1mm]
	\Theta_t = \Theta_{xx} + \gamma(\Theta) u_{xt}^2 - f(\Theta) u_{xt},
	\ear
\ee
which arises in the modeling of one-dimensional thermoviscoelastic processes, where $u=u(x,t)$ 
represents the displacement variable, and $\Theta=\Theta(x,t)$ denotes the temperature.\abs
In fact, in the most prototypical case when $\gamma\equiv const.$ and $f(\Theta)=\Theta$, (\ref{00}) reduces
to the apparently simplest version of the classical model for 
the transfer of mechanical energy into heat due to viscosity effects in solid materials of Kelvin-Voigt type
(\cite{roubicek}, \cite{racke_zheng}).
In some contrast to higher-dimensional counterparts for which the mathematical theory has considerably developed but yet
relies on restrictions on sufficiently sublinear growth of $f$ or also on the presence of additional
dissipative effects (\cite{roubicek}, \cite{mielke_roubicek}, \cite{blanchard_guibe},
\cite{rossi_roubicek_interfaces13}, \cite{gawinecki_zajaczkowski_cpaa}, 
\cite{gawinecki_zajaczkowski}, \cite{roubicek_nodea2013},
\cite{pawlow_zajaczkowski_cpaa17}, \cite{owczarek_wielgos}),
the one-dimensional problem (\ref{00}) as well as some close relatives involving constant viscosities 
have rather comprehensively been investigated with regard both to basic existence theories  
(\cite{dafermos}, \cite{dafermos_hsiao_smooth}, \cite{racke_zheng}, \cite{guo_zhu}, \cite{shen_zheng_zhu},
\cite{kim}, \cite{zheng_shen}, \cite{chen_hoffmann})
and to qualitative behavior (\cite{racke_zheng}, \cite{hsiao_luo}).
Beyond this, results on global existence of solutions at various levels of regularity have also been established
within some classes of viscosities $\gamma$ which are allowed to depend on the strain $u_x$ (\cite{jiang_QAM1993}, \cite{watson})
or also on its the time derivative $u_{xt}$ (\cite{roubicek}).
In the limit case $\gamma\equiv 0$, corresponding problems from thermoelaticity have 
been well-understood not only in one-dimensional domains
(\cite{slemrod}, \cite{bies_cieslak}, \cite{jiang1990}, \cite{racke90}, \cite{racke_shibata},
\cite{racke_shibata_zheng}), but meanwhile also in higher-dimensional situations (\cite{cieslak}).\abs
The knowledge appears significantly sparser, however, in cases when the mechanical processes modeled by (\ref{00})
are considered to depend explicitly on the temperture variable, although both key ingredients 
$\gamma$ and also $a$ have been found to depend on $\Theta$ in various experimental settings of relevance
(\cite{friesen}, \cite{Gubinyi2007}).
In the simultaneous presence of temperature dependencies both in viscosity and in the elastic parameter,
namely, only in certain small-data settings some associated initial-boundary value problems have so far been
proved to admit global solutions (\cite{fricke}, \cite{claes_win}, \cite{win_AMOP}), or at least to possess
solutions defined up to a prescribed finite time (\cite{meyer}); a higher-dimensional version could even be addressed only
under the additional assumption that $f\equiv 0$ (\cite{claes_lankeit_win}). \abs
After all, recent analysis has indicated that a more complete picture can well be expected 
in the case when at least $a\equiv const.$, which seems of particular
relevance for the modeling of viscoelastic dynamics interacting with electric fields in piezoceramic materials (\cite{fricke}, 
\cite{win_global_weak}). 
Indeed, in \cite{win_global_weak} 
it has been seen that in such settings for all temperature-dependent continuous viscosities satisfying 
\be{01}
	0<\inf_{\Theta\ge 0} \gamma(\Theta)<\sup_{\Theta\ge 0} \gamma(\Theta)<\infty,
\ee
certain global solutions to the initial-boundary value problem
\be{0}
	\lball
	u_{tt} = \big(\gamma(\Theta) u_{xt}\big)_x + au_{xx} - \big(f(\Theta)\big)_x,
	\qquad & x\in\Om, \ t>0, \\[1mm]
	\Theta_t = \Theta_{xx} + \gamma(\Theta) u_{xt}^2 - f(\Theta) u_{xt},
	\qquad & x\in\Om, \ t>0, \\[1mm]
	u=0, \quad \Theta_x=0,
	\qquad & x\in\pO, \ t>0, \\[1mm]
	u(x,0)=u_0(x), \quad u_t(x,0)=u_{0t}(x), \quad \Theta(x,0)=\Theta_0(x),
	\qquad & x\in\Om,
	\ear
\ee
posed in the interval $\Om=(0,L)\subset\R$ with $L>0$,
can be constructed for arbitrarily large initial data $u_0, u_{0t}$ and $\Theta_0\ge 0$.
The knowledge on regularity of these solutions, however, is fairly limited and essentially restricted to the information that
\be{02}
	\lbal
	u \in C^0(\bom\times [0,\infty)) \cap L^\infty((0,\infty);W_0^{1,2}(\Om)), \\[1mm]
	u_t \in L^\infty((0,\infty);L^2(\Om)) \cap L^2_{loc}([0,\infty);W_0^{1,2}(\Om))
	\qquad \mbox{and} \\[1mm]
	\Theta\in L^\infty((0,\infty);L^1(\Om)) \cap \bigcap_{q\in [1,3)} L^q_{loc}(\bom\times [0,\infty))
		\cap \bigcap_{r\in [1,\frac{3}{2})} L^r_{loc}([0,\infty);W^{1,r}(\Om)),
	\ear
\ee
which particularly is insufficient to rule out the formation of singular temperature hotspots in the sense of
finite-time blow-up with respect to spatial $L^\infty$ norms in the variable $\Theta$.\abs
{\bf Main results.} \quad
The present manuscript now intends to develop an approach capable of constructing also some regular large-data solutions
to the problem (\ref{0}) with temperature-dependent $\gamma$ at a temporally global level.
At its core, this method will be based on the ambition to further develop some basic {\em a priori} regularity information,
as expressed in the energy conservation law
\be{en}
	\frac{d}{dt} \bigg\{
	\frac{1}{2} \io u_t^2
	+ \frac{a}{2} \io u_x^2
	+ \io \Theta \bigg\}
	= 0,
\ee
by appropriately exploiting the parabolic character of the problem
\be{0z}
	\lball
	z_t = \gamma_0 z_{xx} - \kappa z + h(x,t),
	\qquad & x\in\Om, \ t>0, \\[1mm]
	z(0,t)=z_x(L,t)=0,
	\qquad & t>0,
	\ear
\ee
formally satisfied by 
\be{z}
	z(x,t):=e^{-\kappa t} \int_0^t v(y,t) dy,
	\qquad v:=u_t,
\ee
and
\be{h}
	h(x,t):=\big(\gamma(\Theta)-\gamma_0\big) z_{xx} + ae^{-\kappa t} u_x - e^{-\kappa t} f(\Theta)
\ee
for arbitrary choices of $\gamma_0>0$ and $\kappa\in\R$.
In fact, on the basis of an $L^1$ bound for $\Theta$ implied by (\ref{en}), a Gagliardo-Nirenberg interpolation can be used to 
control the growth of higher-order Lebesgue norms of $\Theta$ in the such a way that for suitable $q>1$, an inequality of the form
\be{03}
	\frac{d}{dt} \io (\Theta+1)^q 
	+ \frac{1}{\Gamma_1} \io (\Theta+1)^{q+2} 
	\le \eta \io v_x^4 + \Gamma_2\|\gamma\|_{L^\infty((0,\infty))}^\lam
\ee
holds for all $\eta>0$, some $\lam=\lam(q)>0$, and some 
$\Gamma_1>0$ and $\Gamma_2>0$ possibly depending on $\eta$, $q$, the initial data and $f$ but
not on $\gamma$, provided that $f$ satisfies an assumption on asymptotic growth that is consistent with 
the prototypical choice $f(\Theta)=\Theta$ (Section \ref{sect3}).
The integral on the right of (\ref{03}), however, can be estimated by means of an argument based on a maximal
Sobolev regularity property of (\ref{0z}): Under an assumption on uniform smallness of the quantity $\gamma-\gamma_0$
entering the crucial second-order contribution to (\ref{h}), namely, for adequately large $\kappa>0$
an appropriately arranged loop-type argument yields bounds, inter alia, for $\int \int v_x^4$ (Section \ref{sect4}).\abs
Using that these, in turn, entail further regularity properties sufficient to ensure global extensibility of local-in-time
solutions available due to known approaches (Lemma \ref{lem_loc}), we can thereby establish the first of our main results
which asserts global classical solvability in (\ref{0}) under mild hypotheses on $f$ and an appropriate smallness
assumption on the total oscillation of $\gamma$, well-consistent with the experimental findings on moderate
temperature dependencies reported in \cite{friesen}:
\begin{theo}\label{theo8}
  Let $L>0$ and $\Om=(0,L) \subset \R$, and let $\gamma_0>0$.
  Then there exists $\del=\del(\gamma_0)>0$ with the property that if $a>0$, if
  \be{gf}
	\lbal
	\gamma \in C^2([0,\infty))
	\qquad \mbox{and} \qquad \\[1mm]
	f\in C^2([0,\infty)) 
	\quad \mbox{satisfies } f(0)=0,
	\ear
  \ee
  and if furthermore
  \be{gg}
	\gamma_0 \le \gamma(\xi) \le \gamma_0+\del
	\qquad \mbox{for all } \xi\ge 0
  \ee
  and
  \be{f}
	|f(\xi)| \le K_f (\xi+1)^\al
	\quad \mbox{for all } \xi\ge 0	
	\quad \mbox{with some $K_f>0$ and $\al\in (0,\frac{3}{2})$,}
  \ee
  then whenever
  \be{init}
	\lbal
	u_0\in W^{3,2}(\Om) \cap W_0^{1,2}(\Om), \\[1mm]
	u_{0t}\in W^{2,2}(\Om) \cap W_0^{1,2}(\Om) \qquad \mbox{and} \\[1mm]
	\Theta_0\in W^{2,2}(\Om)
	\mbox{ satisfies $\Theta_0\ge 0$ in $\Om$ and $\Theta_{0x}=0$ on $\pO$,}
	\ear
  \ee
  there exist uniquely determined functions
  \be{8.1}
	\lbal
	u\in \Big( \bigcup_{\beta\in (0,1)} C^{1+\beta,\frac{1+\beta}{2}}(\bom\times [0,\infty))\Big) 
		\cap C^{2,1}(\bom\times (0,\infty))
		\qquad \mbox{and} \\[1mm]
	\Theta\in \Big( \bigcup_{\beta\in (0,1)} C^{1+\beta,\frac{1+\beta}{2}}(\bom\times [0,\infty))\Big) 
		\cap C^{2,1}(\bom\times (0,\infty))
	\ear
  \ee
  fulfilling
  \be{8.2}
	\begin{array}{l}
	u_t\in 
	\Big( \bigcup_{\beta\in (0,1)} C^{1+\beta,\frac{1+\beta}{2}}(\bom\times [0,\infty))\Big) 
		\cap C^{2,1}(\bom\times (0,\infty)),
	\end{array}
  \ee
  which are such that $\Theta\ge 0$ in $\Om\times (0,\infty)$, and that (\ref{0}) is solved in the classical sense.
\end{theo}
The above analysis can be extended so as to warrant solvability also under less restrictive assumptions on the initial data.
To substantiate this, for conveience in notation we abbreviate
$W^{2,2}_N(\Om):=\{\vp\in W^{2,2}(\Om) \ | \ \vp_x=0 \mbox{ on } \pO\}$,
and specify a notion of solvability that operates at levels of regularity below those from Theorem \ref{theo8}, but
yet markedly above the ones contained in (\ref{02}).
\begin{defi}\label{dw}
  A pair of funtions
  \be{w1}
	\lbal
	u\in C^0(\bom\times [0,\infty)) \cap L^\infty_{loc}([0,\infty);W^{2,2}(\Om)\cap W_0^{1,2}(\Om))
	\qquad \mbox{and} \\[1mm]
	\Theta\in C^0(\bom\times [0,\infty)) \cap L^2_{loc}([0,\infty);W^{2,2}_N(\Om))
	\ear
  \ee
  such that
  \be{w2}
	u_t \in C^0(\bom\times [0,\infty)) \cap L^2_{loc}([0,\infty);W^{2,2}(\Om))
	\qquad \mbox{and} \qquad
	u_{tt} \in L^2_{loc}(\bom\times [0,\infty)),
  \ee
  and that
  \be{w3}
	\Theta_t \in L^2_{loc}(\bom\times [0,\infty)),
  \ee
  will be called a {\em global strong solution} of (\ref{0}) if
  \be{w4}
	u(\cdot,0)=u_0,
	\quad
	u_t(\cdot,0)=u_{0t}
	\quad \mbox{and} \quad
	\Theta(\cdot,0)=\Theta_0
	\qquad \mbox{in } \Om,
  \ee
  and if the identities
  \be{wu}
	u_{tt} = \gamma(\Theta) u_{xxt} + \gamma'(\Theta) \Theta_x u_{xt} + a u_{xx} - f'(\Theta)\Theta_x
  \ee
  and
  \be{wt}
	\Theta_t = \Theta_{xx} + \gamma(\Theta) u_{xt}^2 - f(\Theta) u_{xt}
  \ee
  hold a.e.~in $\Om\times (0,\infty)$.
\end{defi}
Within this framework, actually as a by-product of our analysis leading to Theorem \ref{theo8} we finally obtain 
in Section \ref{sect6} and Section \ref{sect7}
the followig statement on global existence of solutions which, despite reduced assumptions on the initial data, are regular enough
so as to exclude any $L^\infty$ blow-up in finite time:
\begin{theo}\label{theo15}
  Let $\Om=(0,L) \subset \R$ with some $L>0$, let $\gamma_0>0$, and suppose that $\gamma$ and $f$
  satisfy (\ref{gf}), (\ref{gg}) and (\ref{f}) with $\del(\gamma_0)$ as provided by Theorem \ref{theo8}.
  Then for any choice of
  \be{Init}
	\lbal
	u_0\in W^{2,2}(\Om) \cap W_0^{1,2}(\Om), \\[1mm]
	u_{0t}\in W_0^{1,4}(\Om) \qquad \mbox{and} \\[1mm]
	\Theta_0\in W^{1,2}(\Om)
	\mbox{ such that $\Theta_0\ge 0$ in $\Om$,}
	\ear
  \ee
  one can find
  \bas
	\lbal
	u\in C^0(\bom\times [0,\infty)) \cap L^\infty_{loc}((0,\infty);W^{2,2}(\Om) \cap W_0^{1,2}(\Om))
	\qquad \mbox{and} \qquad \\[1mm]
	\Theta\in C^0(\bom\times [0,\infty)) \cap L^\infty_{loc}([0,\infty);W^{1,2}(\Om)) \cap L^2_{loc}([0,\infty);W_N^{2,2}(\Om))
	\ear
  \eas
  such that
  \bas
	u_t\in C^0(\bom\times [0,\infty)) \cap L^\infty_{loc}((0,\infty);W_0^{1,2}(\Om)) 
		\cap L^2_{loc}([0,\infty);W^{2,2}(\Om)),
  \eas
  that $\Theta\ge 0$ in $\Om\times [0,\infty)$, and that $(u,\Theta)$ forms a global strong solution of
  (\ref{0}) in the sense of Definition \ref{dw}.
\end{theo}
\mysection{Preliminaries}\label{sect2}
To begin with, let us record that in the presence of initial data regular enough to comply with
the requirements in Theorem \ref{theo8}, local-in-time classical solutions satisfying a handy extensibility
criterion can always be found.
\begin{lem}\label{lem_loc}
  Let $a>0$, assume that $\gamma$ and $f$ satisfy (\ref{gf}) with $\gamma>0$ on $[0,\infty)$,
  and suppose that (\ref{init}) holds.
  Then there exist $\tm\in (0,\infty]$ as well as uniquely determined functions
  \be{l1}
	\lbal
	u\in \Big( \bigcup_{\beta\in (0,1)} C^{1+\beta,\frac{1+\beta}{2}}(\bom\times [0,\tm))\Big) \cap C^{2,1}(\bom\times (0,\tm))
		\qquad \mbox{and} \\[1mm]
	\Theta\in \Big( \bigcup_{\beta\in (0,1)} C^{1+\beta,\frac{1+\beta}{2}}(\bom\times [0,\tm))\Big) 
		\cap C^{2,1}(\bom\times (0,\tm))
	\ear
  \ee
  which are such that
  \be{l2}
	\begin{array}{l}
	u_t\in 
	\Big( \bigcup_{\beta\in (0,1)} C^{1+\beta,\frac{1+\beta}{2}}(\bom\times [0,\tm))\Big) \cap C^{2,1}(\bom\times (0,\tm)),
	\end{array}
  \ee
  that $\Theta\ge 0$ in $\Om\times (0,\tm)$, that $(u,\Theta)$ solves (\ref{0}) in the classical sense
  in $\Om\times (0,\tm)$, and which additionally are such that
  \bea{ext}
	\mbox{if $\tm<\infty$, \quad then \quad} 
	\limsup_{t\nearrow\tm} \|\Theta(\cdot,t)\|_{W^{1,2}(\Om)} = \infty.
  \eea
\end{lem}
\proof
  This can be obtained by minor adaptation of the reasonings in \cite{fricke} and \cite{win_global_weak}, 
  where statements of this form have been derived for variants of (\ref{0}) in which the present boundary
  condition $u|_{\pO}=0$ has been replaced with the homogeneous Neumann condition $u_x|_{\pO}=0$.
\qed
Throughout the sequel, whenever $a,\gamma,f$ and $(u_0,u_{0t},\Theta_0)$ have been fixed in such a way that the assumptions
of Lemma \ref{lem_loc} are met, we let $\tm$ as well as $u$ and $\Theta$ be as found above, and note
that in this case, the triple $(v,u,\Theta)$ with $v:=u_t$  is a classical solution of
\be{0v}
	\lball
	v_t = \big( \gamma(\Theta) v_x\big)_x + au_{xx} - \big( f(\Theta)\big)_x,
	\qquad & x\in\Om, \ t\in (0,\tm), \\[1mm]
	u_t = v,
	\qquad & x\in\Om, \ t\in (0,\tm), \\[1mm]
	\Theta_t = \Theta_{xx} + \gamma(\Theta) v_x^2 - f(\Theta) v_x,
	\qquad & x\in\Om, \ t\in (0,\tm), \\[1mm]
	v=0, \quad u=0, \quad \Theta_x=0,
	\qquad & x\in\pO, \ t\in (0,\tm), \\[1mm]
	v(x,0)=u_{0t}(x), \quad u(x,0)=u_0(x), \quad \Theta(x,0)=\Theta_0(x),
	\qquad & x\in\Om.
	\ear
\ee
In order to simultaneously address our objectives in Theorem \ref{theo8} and Theorem \ref{theo15},
let us supplement our assumptions in (\ref{init}) by the quantitative requirement that
\be{iM}
	\|u_0\|_{L^\infty(\Om)} 
	+ \|u_{0x}\|_{L^\infty(\Om)}
	+ \|u_{0xx}\|_{L^2(\Om)}
	+ \|u_{0t}\|_{L^\infty(\Om)}
	+ \big\| \big(u_{0t}\big)_x\big\|_{L^4(\Om)}
	+ \|\Theta_0\|_{L^\infty(\Om)}
	+ \|\Theta_{0x}\|_{L^2(\Om)}
	\le M
\ee
with some $M>0$.\abs
In dependence on these bounds, the following information results from the conservation property in (\ref{en}).
\begin{lem}\label{lem1}
  For all $M>0$ there exists $\Lambda_1=\Lambda_1(M)>0$ such that if $a>0$ and (\ref{gf}), (\ref{init}) as well as (\ref{iM})
  hold, then
  \be{1.1}
	\io v^2(\cdot,t)
	\le \Lambda_1
	\qquad \mbox{for all } t\in (0,\tm)
  \ee
  and
  \be{1.2}
	\io u_x^2(\cdot,t) \le \Lambda_1
	\qquad \mbox{for all } t\in (0,\tm)
  \ee
  as well as
  \be{1.3}
	\io \Theta(\cdot,t) \le \Lambda_1
	\qquad \mbox{for all } t\in (0,\tm).
  \ee
\end{lem}
\proof
  As can readily be verified on the basis of (\ref{0v}) and integration by parts,
  \bas
	\frac{d}{dt} \io v^2
	&=& - \io \gamma(\Theta) v_x^2 - a \io u_x v_x
	+ \io f(\Theta) v_x \\
	&=& - \io \gamma(\Theta) v_x^2 - \frac{a}{2} \frac{d}{dt} \io u_x^2
	+ \io f(\Theta) v_x 
  \eas
  and
  \bas
	\frac{d}{dt} \io \Theta = \io \gamma(\Theta) v_x^2 - \io f(\Theta) v_x
  \eas
  for all $t\in (0,\tm)$.
  On adding these identities and integrating in time, we see that
  \bas
	\frac{1}{2} \io v^2
	+ \frac{a}{2} \io u_x^2
	+ \io \Theta
	= \frac{1}{2} \io u_{0t}^2
	+ \frac{a}{2} \io u_{0x}^2
	+ \io \Theta_0
	\qquad \mbox{for all } t\in (0,\tm),
  \eas
  whence using (\ref{iM}) in estimating
  \bas
	\frac{1}{2} \io u_{0t}^2
	+ \frac{a}{2} \io u_{0x}^2
	+ \io \Theta_0
	\le
	\frac{M^2|\Om|}{2} + \frac{aM^2 |\Om|}{2} + M|\Om|
  \eas
  we obtain (\ref{1.1})-(\ref{1.3}).	
\qed
\mysection{Controlling the growth of $\io (\Theta+1)^q$ in terms of $\io v_x^4$}\label{sect3}
In order to derive regularity properties beyond those from Lemma \ref{lem1}, let us first apply a standard
variational argument to the third equation in (\ref{0v}).
By means of an appropriate interpolation relying on (\ref{1.3}) at its lower end, in a first step
in this regard we obtain the following yet fairly general statement on control of $\Theta+1$ in dependence
on certain Lebesgue norms of the expressions $v_x$ which act as an essential part of the inhomogeneity in the
temperature evolution modeled by (\ref{0v}).
\begin{lem}\label{lem21}
  Let $f\in C^2([0,\infty))$ satisfy (\ref{f}), and let $M>0$ and $q>1$.
  Then there exists $\Gamma_1=\Gamma_1(M,q,f)>0$ such that 
  whenever $a>0$, $0<\gamma\in C^2([0,\infty)) \cap L^\infty((0,\infty))$ and (\ref{init}) as well as (\ref{iM}) hold,
  \bea{21.1}
	\frac{d}{dt} \io (\Theta+1)^q
	+ \frac{1}{\Gamma_1} \io (\Theta+1)^{q+2}
	&\le& \Gamma_1 \cdot \|\gamma\|_{L^\infty((0,\infty))}^\frac{q+2}{3} \io |v_x|^\frac{2(q+2)}{3} \nn\\
	& & + \Gamma_1 \io |v_x|^\frac{q+2}{3-\al}
	+ \Gamma_1
	\qquad \mbox{for all } t\in (0,\tm).
  \eea
\end{lem}
\proof
  According to the Gagliardo-Nirenberg inequality, there exists $c_1=c_1(q)>0$ such that
  \be{21.2}
	\|\vp\|_{L^\frac{2(q+2)}{q}(\Om)}^\frac{2(q+2)}{q} 
	\le c_1 \|\vp_x\|_{L^2(\Om)}^2 \|\vp\|_{L^\frac{2}{q}(\Om)}^\frac{4}{q}
	+ c_1 \|\vp\|_{L^\frac{2}{q}(\Om)}^\frac{2(q+2)}{q}
	\qquad \mbox{for all } \vp\in W^{1,2}(\Om),
  \ee
  and abbreviating
  \be{21.3}
	c_2\equiv c_2(M,q):=\frac{q^2 c_1}{4} \cdot \big( |\Om| + \Lambda_1\big)^2
	\qquad \mbox{and} \qquad
	c_3\equiv c_3(M,q):=c_1 \cdot \big( |\Om| + \Lambda_1 \big)^{q+2}
  \ee
  as well as
  \be{21.4}
	c_3\equiv c_3(M,q):=\frac{q(q-1)}{c_2}
	\qquad \mbox{and} \qquad
	c_4\equiv c_4(M,q):=\frac{q(q-1)c_3}{c_2},
  \ee
  with $\Lambda_1=\Lambda_1(M)$ taken from Lemma \ref{lem1}, we claim that the intended conclusion holds if we let
  \be{21.5}
	\Gamma_1\equiv \Gamma_1(M,q,f)
	:=\max \bigg\{ \frac{2}{c_3} \, , \, \Big(\frac{4}{c_3}\Big)^\frac{q-1}{3} \cdot q^\frac{q+2}{3} \, , \,
		\Big(\frac{4}{c_3}\Big)^\frac{q+\al-1}{3-\al} \cdot (qK_f)^\frac{q+2}{3-\al} \, , \, c_4 \bigg\}.
  \ee
  To see this, given $a>0$ and $0<\gamma\in C^2([0,\infty))\cap L^\infty((0,\infty))$ as well as initial data fulfilling
  (\ref{init}) and (\ref{iM}), we use the third equation in (\ref{0v})
  and integrate by parts to compute
  \bea{21.6}
	\frac{d}{dt} \io (\Theta+1)^q
	&=& q\io (\Theta+1)^{q-1} \cdot \big\{ \Theta_{xx} + \gamma(\Theta) v_x^2 - f(\Theta) v_x \big\} \nn\\
	&=& - q(q-1) \io (\Theta+1)^{q-2} \Theta_x^2
	+ q \io \gamma(\Theta) (\Theta+1)^{q-1} v_x^2 \nn\\
	& & - q \io f(\Theta) (\Theta+1)^{q-1} v_x
	\qquad \mbox{for all } t\in (0,\tm),
  \eea
  we we note that due to (\ref{21.2}), (\ref{1.1}) and (\ref{21.3}),
  \bas
	\io (\Theta+1)^{q+2}
	&=& \big\| (\Theta+1)^\frac{q}{2} \big\|_{L^\frac{2(q+2)}{q}(\Om)}^\frac{2(q+2)}{q} \\
	&\le& c_1 \big\| \big( (\Theta+1)^\frac{q}{2}\big)_x \big\|_{L^2(\Om)}^2
		\big\| (\Theta+1)^\frac{q}{2} \big\|_{L^\frac{2}{q}(\Om)}^\frac{4}{q}
	+ c_1\big\| (\Theta+1)^\frac{2}{2} \big\|_{L^\frac{2}{q}(\Om)}^\frac{2(q+2)}{q} \\
	&=& \frac{q^2 c_1}{4} \cdot \bigg\{ \io (\Theta+1)^{q-2} \Theta_x^2 \bigg\} 
		\cdot \bigg\{ \io (\Theta+1) \bigg\}^2
	+ c_1 \cdot\bigg\{ \io (\Theta+1) \bigg\}^{q+2} \\
	&\le& \frac{q^2 c_1}{4} \cdot \bigg\{ \io (\Theta+1)^{q-2} \Theta_x^2 \bigg\}		
		\cdot \big( |\Om|+ \Lambda_1\big)^2
	+ c_1\cdot \big( |\Om| + \Lambda_1\big)^{q+2} \\
	&=& c_2 \io (\Theta+1)^{q-2} \Theta_x^2
	+ c_3
	\qquad \mbox{for all } t\in (0,\tm).
  \eas
  As thus, in line with (\ref{21.4}),
  \bas
	q(q-1) \io (\Theta+1)^{q-2} \Theta_x^2
	&\ge& q(q-1) \cdot \bigg\{ \frac{1}{c_2} \io (\Theta+1)^{q+2} - \frac{c_3}{c_2}\bigg\} \\
	&=& c_3 \io (\Theta+1)^{q+2} - c_4
	\qquad \mbox{for all } t\in (0,\tm),
  \eas
  from (\ref{21.6}) we obtain that
  for all $t\in (0,\tm)$,
  \be{21.7}
	\frac{d}{dt} \io (\Theta+1)^q
	+ c_3 \io (\Theta+1)^{q+2}
	\le c_4 + q\io \gamma(\Theta) (\Theta+1)^{q-1} v_x^2
	- q \io f(\Theta) (\Theta+1)^{q-1} v_x.
  \ee
  Here, Young's inequality guarantees that
  \bas
	q \io \gamma(\Theta) (\Theta+1)^{q-1} v_x^2
	&\le& q\|\gamma\|_{L^\infty((0,\infty))} \io (\Theta+1)^{q-1} v_x^2 \\
	&=& \io \Big\{ \frac{c_3}{4} (\Theta+1)^{q+2} \Big\}^\frac{q-1}{q+2} 
		\cdot \Big(\frac{4}{c_3}\Big)^\frac{q-1}{q+2} \cdot q \|\gamma\|_{L^\infty((0,\infty))} v_x^2 \\
	&\le& \frac{c_3}{4} \io (\Theta+1)^{q+2}  \\
	& & + \Big(\frac{4}{c_3}\Big)^\frac{q-1}{3} \cdot q^\frac{q+2}{3} \|\gamma\|_{L^\infty((0,\infty))}^\frac{q+2}{3}
		\io |v_x|^\frac{2(q+2)}{3}
	\qquad \mbox{for all } t\in (0,\tm),
  \eas
  and that by (\ref{f}), as we are particularly assuming that $\al<3$,
  \bas
	- q \io f(\Theta) (\Theta+1)^{q-1} v_x
	&\le& q K_f \io (\Theta+1)^{q+\al-1} |v_x| \\
	&=& \io \Big\{ \frac{c_3}{4} (\Theta+1)^{q+2} \Big\}^\frac{q+\al-1}{q+2} 	
		\cdot \Big(\frac{4}{c_3}\Big)^\frac{q+\al-1}{q+2} \cdot q K_f |v_x| \\
	&\le& \frac{c_3}{4} \io (\Theta+1)^{q+2} \\
	& & + \Big(\frac{4}{c_3}\Big)^\frac{q+\al-1}{3-\al} \cdot (q K_f)^\frac{q+2}{3-\al} 
		\cdot \io |v_x|^\frac{q+2}{3-\al}
	\qquad \mbox{for all } t\in (0,\tm).
  \eas
  Therefore, (\ref{21.7}) entails that
  \bas
	\frac{d}{dt} \io (\Theta+1)^q
	+ \frac{c_3}{2} \io (\Theta+1)^{q+2}
	&\le& c_4 
	+ \Big(\frac{4}{c_3}\Big)^\frac{q-1}{3} \cdot q^\frac{q+2}{3} \|\gamma\|_{L^\infty((0,\infty))}^\frac{q+2}{3} 
		\io |v_x|^\frac{2(q+2)}{3} \nn\\
	& & + \Big(\frac{4}{c_3}\Big)^\frac{q+\al-1}{3-\al} \cdot (q K_f)^\frac{q+2}{3-\al} \io |v_x|^\frac{q+2}{3-\al}
	\qquad \mbox{for all } t\in (0,\tm),
  \eas
  which in view of (\ref{21.5}) implies (\ref{21.1}).
\qed
Within a certain range of superlinear powers $q$, one further and quite simple interpolation shows that both integrals on the right
of (\ref{21.1}) can be bounded against, essentially, a small portion of the quantity $\io v_x^4$ that will later on be
controlled in the course of a loop-type argument involving a maximal Sobolev regularity estimate (see Lemma \ref{lem4}).
\begin{lem}\label{lem2}
  Let $f\in C^2([0,\infty))$ be such that (\ref{f}) holds, and let $M>0$.
  Then whenever $q>1$ satisfies
  \be{2.1}
	q<4
	\qquad \mbox{and} \qquad
	q<10-4\al,
  \ee
  given any $\eta>0$ one can choose $\Gamma_2=\Gamma_2(M,\eta,q,f)>0$ in such a way that if $a>0$ and
  $0<\gamma\in C^2([0,\infty))\cap L^\infty((0,\infty))$, and if (\ref{init}) and (\ref{iM}) hold, it follows that
  \be{2.3}
	\frac{d}{dt} \io (\Theta+1)^q 
	+ \frac{1}{\Gamma_1} \io (\Theta+1)^{q+2} 
	\le \eta \io v_x^4 + \Gamma_2\|\gamma\|_{L^\infty((0,\infty))}^\frac{2(q+2)}{4-q}
	+ \Gamma_2
	\qquad \mbox{for all } t\in (0,\tm),
  \ee
  where $\Gamma_1=\Gamma_1(M,q,f)$ is as in Lemma \ref{lem21}.
\end{lem}
\proof
  We let $\Gamma_1=\Gamma_1(M,q,f)$ be as in Lemma \ref{lem21}, and noting that $4-q$ and $10-q-4\al$ are both 
  positive by (\ref{2.1}), given $\eta>0$ we let
  \be{2.4}
	\Gamma_2\equiv \Gamma_2(M,\eta,q,f)
	:= \max \bigg\{ \Gamma_1 \, , \, \Big(\frac{2}{\eta}\Big)^\frac{q+2}{4-q} \Gamma_1^\frac{6}{4-q} |\Om| \, , \,
		\Big(\frac{2}{\eta}\Big)^\frac{q+2}{10-q-4\al} \Gamma_1^\frac{12-4\al}{10-q-4\al} |\Om| \bigg\}.
  \ee
  Then fixing $a>0$ and $0<\gamma\in C^2([0,\infty)) \cap L^\infty((0,\infty))$, and assuming that (\ref{init})
  and (\ref{iM}) hold, in (\ref{21.1}) we can rely on (\ref{2.1}) to see by means of Young's inequality that
  \bas
	\Gamma_1 \cdot \|\gamma\|_{L^\infty((0,\infty))}^\frac{q+2}{3} \io |v_x|^\frac{2(q+2)}{3}
	&=& \io \Big\{ \frac{\eta}{2} v_x^4 \Big\}^\frac{q+2}{6} 
		\cdot \Big(\frac{2}{\eta}\Big)^\frac{q+2}{6} \Gamma_1 \|\gamma\|_{L^\infty((0,\infty))}^\frac{q+2}{3} \\
	&\le& \frac{\eta}{2} \io v_x^4
	+ \Big(\frac{2}{\eta}\Big)^\frac{q+2}{4-q} \Gamma_1^\frac{6}{4-q} |\Om| \|\gamma\|_{L^\infty((0,\infty))}^\frac{2(q+2)}{4-q}
	\qquad \mbox{for all } t\in (0,\tm),
  \eas
  and that, similarly,
  \bas
	\Gamma_1 \io |v_x|^\frac{q+2}{3-\al}
	&=& \io \Big\{ \frac{\eta}{2} v_x^4 \Big\}^\frac{q+2}{12-4\al} 
		\cdot \Big(\frac{2}{\eta}\Big)^\frac{q+2}{12-4\al} \Gamma_1 \\
	&\le& \frac{\eta}{2} \io v_x^4
	+ \Big(\frac{2}{\eta}\Big)^\frac{q+2}{10-q-4\al} \Gamma_1^\frac{12-4\al}{10-q-4\al} |\Om|
	\qquad \mbox{for all } t\in (0,\tm).
  \eas
  In line with (\ref{2.4}), from (\ref{21.1}) we thus obtain (\ref{2.3}).
\qed
\mysection{Space-time bounds for $v_x^4$ and $\Theta^{4\al}$ via maximal Sobolev regularity}\label{sect4}
The key part of our analysis now relies on the parabolic character of the first equation in (\ref{0v})
to relate the expression on the right-hand side of (\ref{2.3}) to the regularity of sources in the rearranged
version (\ref{z})-(\ref{0z}) of this evolution problem.\abs
Our considerations in this direction are based on the following consequence of general theory on maximal Sobolev
regularity for the scalar heat equation, augmented by a simple reflection argument that provides accessibility of
the mixed boundary value problem in (\ref{0z}) to classical statements addressing homogeneous Dirichlet data.
\begin{lem}\label{lem3}
  For each $p\in (1,\infty)$ and any $D>0$ there exists $K(p,D)>0$ with the property that whenever
  $T>0, g\in C^0(\bom\times [0,T])$ and $z\in C^{2,1}(\bom\times [0,T])$ are such that
  \be{3.1}
	\lball
	z_t = D z_{xx} + g(x,t),
	\qquad & x\in\Om, \ t\in (0,T), \\[1mm]
	z(0,t)=z_x(L,t)=0,
	\qquad & t\in (0,T),
	\ear
  \ee
  is solved in the classical sense, it follows that
  \be{3.2}
	\int_0^T \io |z_{xx}|^p \le K(p,\gamma_0) \io |z_{xx}(x,0)|^p
	+ K(p,\gamma_0) \int_0^T \io |g|^p.
  \ee
\end{lem}
\proof
  As a consequence of a well-known result on maximal Sobolev regularity in the Dirichlet problem for the inhomogeneous
  heat equation (\cite{giga_sohr}, \cite{hieber_pruess}), one can choose $c_1=c_1(p,D)>0$ in such a way that
  of $T>0, \wh{g}\in C^0([0,2L]\times [0,T])$ and $\wh{z}\in C^{2,1}([0,2L]\times [0,T])$ satisfy
  \bas
	\lball
	\wh{z}_t = D \wh{z}_{xx} + \wh{g}(x,t),
	\qquad & x\in (0,2L), \ t\in (0,T), \\[1mm]
	\wh{z}(0,t)=\wh{z}(L,t)=0,
	\qquad & t\in (0,T),
	\ear
  \eas
  then
  \bas
	\int_0^T \int_0^{2L} |\wh{z}_{xx}|^p	
	\le c_1 \int_0^{2L} |\wh{z}_{xx}(\cdot,0)|^p
	+ c_1 \int_0^T \int_0^{2L} |\wh{g}|^p.
  \eas
  For $T,g$ and $z$ with the listed properties, this yields (\ref{3.2}) with $K(p,\gamma_0):=c_1$, as in this situation
  the above can be applied to the functions $\wh{g}$ and $\wh{z}$ defined via reflection according to
  \bas
	(\wh{g},\wh{z})(x,t):=\lball
	(g,z)(x,t),
	\qquad & x\in [0,L], \\[1mm]
	(g,z)(2L-x,t),
	\qquad & x\in (L,2L], 	
	\ear
  \eas
  for $t\in [0,T]$.
\qed
Now in the range of $\al$ specified in Theorem \ref{theo8}, the core of our reasoning reveals the following
by a combination of Lemma \ref{lem2} with Lemma \ref{lem3},
the latter being applied to $p:=4$ in view of the particular integrability power appearing on the right of (\ref{2.3}).
\begin{lem}\label{lem4}
  For each $\gamma_0>0$ there exists $\del=\del(\gamma_0)>0$ with the property that if $a>0$ 
  and (\ref{gf}), (\ref{gg}) as well as (\ref{f}) are valid,
  given any $M>0$ and $T_0>0$ one can find $\Lambda_2=\Lambda_2(M,T_0,a,\gamma,f)>0$ such that
  if (\ref{init}) and (\ref{iM}) holds, we have
  \be{4.1}
	\int_0^T \io v_x^4 \le \Lambda_2
	\qquad \mbox{for all } T\in (0,T_0) \cap (0,\tm)
  \ee
  as well as
  \be{4.01}
	\int_0^T \io (\Theta+1)^{4\al} \le \Lambda_2
	\qquad \mbox{for all } T\in (0,T_0) \cap (0,\tm).
  \ee
\end{lem}
\proof
  Without loss of generality assuming that $\al>\frac{3}{4}$, we let $K:=K(4,\gamma_0)$ with 
  $(K(p,D))_{p>1,D>0}$ as in Lemma \ref{lem3},
  and choose $\del=\del(\gamma_0)>0$ small enough such that
  \be{4.3}
	256 K \del^4 \le \frac{1}{4}.
  \ee
  Then assuming that $a>0$ and that (\ref{gf}), (\ref{gg}), (\ref{f}) and (\ref{init}) as well as (\ref{iM}) 
  hold with some $M>0$, fixing $T_0>0$
  we choose $\kappa=\kappa(T_0,a)>0$ suitably large such that
  \be{4.4}
	\frac{27\cdot 32 K a^4 T_0}{\kappa^3} \le \frac{1}{4},
  \ee
  and we thereupon pick some small $\eta=\eta(M,T_0,a,f)>0$ fulfilling
  \be{4.6}
	4^4 K K_f^4 \Gamma_1 \eta \le \frac{1}{4e^{4\kappa T_0}},
  \ee
  where $\Gamma_1=\Gamma_1(M,q,f)$ is as in Lemma \ref{lem21}, with
  \be{q}
	q:=4\al-2
  \ee
  satisfying $q>1$ according to our assumption $\al>\frac{3}{4}$.\abs
  Now to derive (\ref{4.1}) and (\ref{4.01}), we let
  \be{4.7}
	V(x,t):=\int_0^x v(y,t) dy,
	\qquad x\in \bom, \ t\in [0,\tm),
  \ee
  and
  \be{4.8}
	z(x,t):= e^{-\kappa t} V(x,t),
	\qquad x\in\bom, \ t\in [0,\tm),
  \ee
  as well as 
  \be{4.9}
	z_0(x):=\int_0^x u_{0t}(y) dy,
	\qquad x\in\bom,
  \ee
  and observe that according to Lemma \ref{lem_loc}, $V$ and hence also $z$ belongs to $C^{2,1}(\bom\times [0,\tm))$
  with $V_x=v$ and $V_{xx}=v_x$, and with $V_t(x,t)=\int_0^x v_t(y,t) dy$ for $x\in\bom$ and $t\in [0,\tm)$, so that by (\ref{0v}),
  \bas
	V_t=\gamma(\Theta) v_x + au_x - f(\Theta)
	= \gamma(\Theta) V_{xx} + au_x - f(\Theta)
	\qquad \mbox{in } \Om\times (0,\tm)
  \eas
  and thus
  \bas
	z_t &=& e^{-\kappa t} \cdot \big\{ \gamma(\Theta) V_{xx} + au_x - f(\Theta)\big\} - \kappa e^{-\kappa t} V \\
	&=& \gamma(\Theta) z_{xx}+ ae^{-\kappa t} u_x
	- e^{-\kappa t} f(\Theta)
	- \kappa z
	\qquad \mbox{in } \Om\times (0,\tm).
  \eas
  Artificially decomposing the diffusive part here, we see that $z$ forms a classical solution of
  \bas
	\lball
	z_t = \gamma_0 z_{xx} + \big(\gamma(\Theta)-\gamma_0\big) z_{xx} + ae^{-\kappa t} u_x
	- e^{-\kappa t} f(\Theta) - \kappa z,
	\qquad & x\in\Om, \ t\in (0,\tm), \\[1mm]
	z(0,t)=z_x(L,t)=0,
	\qquad & t\in (0,\tm),
	\ear
  \eas
  so that Lemma \ref{lem3} applies so as to assert that in line with our definition of $K$,
  \bea{4.10}
	\int_0^T \io z_{xx}^4
	&\le& K \io z_{0xx}^4
	+ K \int_0^T \io \Big| \big(\gamma(\Theta)-\gamma_0\big) z_{xx} + ae^{-\kappa t} u_x - e^{-\kappa t} f(\Theta) - \kappa z
		\Big|^4 \nn\\
	&\le& K\io z_{0xx}^4
	+ 256 K \int_0^T \io \big| \gamma(\Theta)-\gamma_0\big|^4 z_{xx}^4
	+ 256 K a^4 \int_0^T \io e^{-4\kappa t} u_x^4 \nn\\
	& & + 256 K \int_0^T \io e^{-4\kappa t} f^4(\Theta) 
	+ 256 K \kappa^4 \int_0^T \io z^4
	\qquad \mbox{for all } T\in (0,\tm),
  \eea
  because $(\xi_1+\xi_2+\xi_3+\xi_4)^4 \le 256 (\xi_1^4 + \xi_2^4 + \xi_3^4 + \xi_4^4)$ for all $(\xi_1,\xi_2,\xi_3,\xi_4)\in\R^4$.
  Here, (\ref{4.9}) and (\ref{iM}) ensure that
  \be{4.111}
	K \io z_{0xx}^4
	= K \io (u_{0t})_x^4
	\le K M^4,
  \ee
  and since
  \bas
	z^4(x,t) \le V^4(x,t)
	= \bigg| \int_0^x v(y,t) dy \bigg|^4
	\le |\Om|^2 \cdot \bigg\{ \io v^2(y,t) dy \bigg\}^2
	\le |\Om|^2 \Lam_1^2
	\quad \mbox{for all $(x,t)\in \Om\times (0,\tm)$}
  \eas
  by the Cauchy-Schwarz inequality and (\ref{1.1}), we obtain that
  \be{4.112}
	256 K \kappa^4 \int_0^T \io z^4
	\le 256K |\Om|^3 \Lam_1^2 \kappa^4 T_0
	\qquad \mbox{for all } T\in (0,T_0)\cap (0,\tm). 
  \ee
  Moreover, (\ref{gg}) along with (\ref{4.3}) guarantees that
  \bea{4.11}
	256K \int_0^T \io \big| \gamma(\Theta)-\gamma_0\big|^4 z_{xx}^4
	&\le& 256 K \del^4 \int_0^T \io z_{xx}^4 \nn\\
	&\le& \frac{1}{4} \int_0^T \io z_{xx}^4
	\qquad \mbox{for all } T\in (0,\tm),
  \eea
  while using (\ref{f}) together with (\ref{q}) we see that
  \bea{4.12}
	256 K \int_0^T \io e^{-4\kappa t} f^4(\Theta)
	&\le& 256 K \int_0^T \io f^4(\Theta) \nn\\
	&\le& 256K K_f^4 \int_0^T \io (\Theta+1)^{4\al} \nn\\
	&=& 256 K K_f^4 \int_0^T \io (\Theta+1)^{q+2}
	\qquad \mbox{for all } T\in (0,\tm).
  \eea
  Now in order to estimate the crucial third to last summand in (\ref{4.10}), we first observe that in view of (\ref{4.7})
  and (\ref{4.8}),
  \bas
	u_x(x,t)
	&=& u_{0x}(x)
	+ \int_0^t v_x(x,s) ds \nn\\
	&=& u_{0x}(x) + \int_0^t e^{\kappa s} z_{xx}(x,s) ds
	\qquad \mbox{for all $x\in\Om$ and } t\in (0,\tm),
  \eas
  so that due to the H\"older inequality,
  \bas
	u_x^4(x,t)
	&\le& 8 u_{0x}^4(x)
	+ 8 \bigg| \int_0^t e^{\kappa s} z_{xx}(x,s) ds \bigg|^4 \\
	&\le& 8 u_{0x}^4(x)
	+ 8 \cdot \bigg\{ \int_0^t z_{xx}^4(x,s) ds \bigg\} \cdot \bigg\{ \int_0^t e^\frac{4\kappa s}{3} ds \bigg\}^3
	\qquad \mbox{for all $x\in\Om$ and } t\in (0,\tm).
  \eas
  Since
  \bas
	\bigg\{ \int_0^t e^\frac{4\kappa s}{3} ds \bigg\}^3
	&=& \bigg\{ \frac{3}{4\kappa} \cdot \Big( e^\frac{4\kappa t}{3} -1\Big) \bigg\}^3
	\le \Big(\frac{4}{4\kappa}\Big)^3 e^{4\kappa t}
	\qquad \mbox{for all } t>0,
  \eas
  this implies that, by (\ref{iM}),
  \bas
	u_x^4(x,t)
	&\le& 8 u_{0x}^4(x)
	+ 8\cdot\Big(\frac{3}{4\kappa}\Big)^3 e^{4\kappa t} \cdot \int_0^t z_{xx}^4(x,s) dxds \nn\\
	&\le& 8 M^4
	+ 8\cdot\Big(\frac{3}{4\kappa}\Big)^3 e^{4\kappa t} \cdot \int_0^t z_{xx}^4(x,s) dxds
	\qquad \mbox{for all $x\in\Om$ and } t\in (0,\tm),
  \eas
  and that thus
  \bea{4.14}
	& & \hs{-30mm}
	256 K a^4 \int_0^T \io e^{-4\kappa t} u_x^4 \nn\\
	&\le& 2048 K a^4 M^4 |\Om| T_0
	+ \frac{27\cdot 32 K a^4}{\kappa^3} \int_0^T \io \bigg\{ \int_0^T z_{xx}^4(x,s) ds \bigg\} dxdt \nn\\
	&\le& 2048 K a^4 \int_0^T \io u_{0x}^4
	+ \frac{27\cdot 32 K a^4 T_0}{\kappa^3} \int_0^T \io z_{xx}^4 \nn\\
	&\le& 2048 K a^4 \int_0^T \io u_{0x}^4 
	+ \frac{1}{4} \int_0^T \io z_{xx}^4
	\qquad \mbox{for all } T\in (0,T_0)\cap (0,\tm)
  \eea
  thanks to our restriction on $\kappa$ in (\ref{4.4}).\abs
  In summary, from (\ref{4.111})-(\ref{4.14}) we infer that (\ref{4.10}) implies the inequality
  \be{4.15}
	\frac{1}{2} \int_0^T \io z_{xx}^4
	\le 256 K K_f^4 \int_0^T \io (\Theta+1)^{q+2}
	+ c_1
	\qquad \mbox{for all } T\in (0,T_0)\cap (0,\tm),
  \ee
  where
  \bas
	c_1 \equiv c_1(M,T_0,a,f)
	:= K M^4
	+ 256 K |\Om|^3 \Lambda_1^2 \kappa^4 T_0
	+ 2048 K a^4 M^4 |\Om| T_0,
  \eas
  with $\Lambda_1=\Lambda_1(M)$ as in Lemma \ref{lem1}.\abs
  To appropriately compensate the integral on the right of (\ref{4.15}), we next observe that (\ref{q}) together with
  our assumption $\al<\frac{3}{2}$ asserts that
  \bas
	q<4\cdot\frac{3}{2} - 2 =4
	\qquad \mbox{and} \qquad
	q-(10-4\al)=8\al-12 <0,
  \eas
  meaning that (\ref{2.1}) is valid, and that we may thus draw on Lemma \ref{lem2} to see upon integrating (\ref{2.3}) that with
  $\Gamma_2=\Gamma_2(M,\eta,q,f)$ as provided there,
  \bas
	\frac{1}{\Gamma_1} \int_0^T \io (\Theta+1)^{q+2}
	\le \io (\Theta_0+1)^q
	+ \eta \int_0^T \io v_x^4
	+ \Gamma_2 \|\gamma\|_{L^\infty(\Om)}^\frac{2(q+2)}{4-q} \cdot T
	+ \Gamma_2 \cdot T
	\qquad \mbox{for all } T\in (0,\tm).
  \eas
  Therefore, as a consequence of our smallness condition on $\eta$ in (\ref{4.6}) we obtain on writing
  \bas
	c_2 \equiv c_2(M,T_0,a,\gamma,f):=256 K K_f^4 \Gamma_1 \cdot \bigg\{ \io (\Theta_0+1)^q
	+ \Gamma_2 \|\gamma\|_{L^\infty((0,\infty))}^\frac{2(q+2)}{4-q} \cdot T_0 + \Gamma_2 T_0 \bigg\}
  \eas
  that, again by (\ref{4.7}) and (\ref{4.8}),
  \bea{4.16}
	256 K K_f^4 \int_0^T \io (\Theta+1)^{q+2}
	&\le& 256 K K_f^4 \Gamma_1 \eta \int_0^T \io v_x^4
	+ c_2 \nn\\
	&\le& \frac{1}{4 e^{4\kappa T_0}} \int_0^T \io v_x^4
	+ c_2
	\qquad \mbox{for all } T\in (0,T_0)\cap (0,\tm).
  \eea
  Since 
  \bas
	\int_0^T \io z_{xx}^4
	= \int_0^T \io e^{-4\kappa t} v_x^4
	\ge e^{-4\kappa T_0} \int_0^T \io v_x^4
	\qquad \mbox{for all } T\in (0,T_0)\cap (0,\tm).
  \eas
  by (\ref{4.7}) and (\ref{4.8}),
  a combination of (\ref{4.16}) with (\ref{4.15}) shows that
  \bas
	\frac{1}{2} e^{-4\kappa T_0} \int_0^T \io v_4^4
	\le \frac{1}{4} e^{-4\kappa T_0} \int_0^T \io v_4^4 + c_1 + c_2
	\qquad \mbox{for all } T\in (0,T_0)\cap (0,\tm).
  \eas
  and that thus
  \bas
	\int_0^T \io v_x^4 \le 4 e^{4\kappa T_0} (c_1+c_2)
	\qquad \mbox{for all } T\in (0,T_0)\cap (0,\tm).
  \eas
  Once more going back to (\ref{4.16}) and recalling (\ref{q}), from this we obtain both (\ref{4.1}) and (\ref{4.01})
  with some adequately large $\Lambda_2=\Lambda_2(M,T_0,a,\gamma,f)>0$.
\qed
\mysection{Higher regularity properties of $\Theta$. Proof of Theorem \ref{theo8}}\label{sect5}
As the estimates in Lemma \ref{lem4} particularly imply space-time $L^2$ bounds for the sources in the heat equation
governing $\Theta$, the following can be obtained by another standard testing procedure.
\begin{lem}\label{lem5}
  Let $\gamma_0>0$, let $\del(\gamma_0)$ be as in Lemma \ref{lem4}, let $a>0$, and assume (\ref{gf}), (\ref{gg}) and (\ref{f}).
  Then for all $M>0$ and each $T_0>0$ there exists $\Lambda_3=\Lambda_3(M,T_0,a,\gamma,f)>0$ 
  such that whenever (\ref{init}) and (\ref{iM}) hold, we have
  \be{5.1}
	\io \Theta_x^2(x,t) dx \le \Lambda_3
	\qquad \mbox{for all } t\in (0,T_0)\cap (0,\tm)
  \ee
  and
  \be{5.01}
	\int_0^T \io \Theta_{xx}^2 \le \Lambda_3
	\qquad \mbox{for all } T\in (0,T_0)\cap (0,\tm).
  \ee
\end{lem}
\proof
  In the identity
  \be{5.2}
	\Theta_t = \Theta_{xx} + g(x,t),
	\qquad x\in\Om, \ t\in (0,\tm),
  \ee
  we can use Young's inequality along with (\ref{f}) to estimate
  \bas
	g:=\gamma(\Theta) v_x^2 + f(\Theta) v_x
  \eas
  according to
  \bas
	g^2 &\le& 2\cdot\big( \gamma^2(\Theta) v_x^4 + f^2(\Theta) v_x^2 \big) \\
	&\le& 2\big(\|\gamma\|_{L^\infty(\Om)}^2 + 1\big) v_x^4 + f^4(\Theta) \\
	&\le& 2\big(\|\gamma\|_{L^\infty(\Om)}^2 + 1\big) v_x^4 + K_f^4 (\Theta+1)^{4\al}
	\qquad \mbox{in } \Om\times (0,\tm).
  \eas
  Given $M>0$ and $T_0>0$,
  by means of Lemma \ref{lem4} we thus obtain that 
  on taking $\Lambda_2=\Lambda_2(M,T_0,a,\gamma,f)$ as provided there, that 
  whenever (\ref{init}) and (\ref{iM}) hold, we have
  \be{5.3}
	\int_0^T \io g^2 
	\le c_1\equiv c_1(M,T_0,a,\gamma,f) :=
	2\big(\|\gamma\|_{L^\infty(\Om)}^2 + 1\big) \Lambda_2 + K_f^4 \Lambda_2
	\qquad \mbox{for all } T\in (0,T_0)\cap (0,\tm),
  \ee
  from which (\ref{5.1}) results in a straightforward manner:
  Testing (\ref{5.2}) against $-\Theta_{xx}$, namely, we obtain by again using Young's inequality that
  \bas
	\frac{1}{2} \frac{d}{dt} \io \Theta_x^2
	+ \io \Theta_{xx}^2
	= - \io \Theta_{xx}
	\le \frac{1}{2} \io \Theta_{xx}^2 + \frac{1}{2} \io g^2
	\qquad \mbox{for all } t\in (0,\tm),
  \eas
  so that
  \bas
	\frac{d}{dt} \io \Theta_x^2
	+ \io \Theta_{xx}^2
	\le \io g^2
	\qquad \mbox{for all } t\in (0,\tm)
  \eas
  and hence, by (\ref{5.3}) and (\ref{iM}) and the continuity of $[0,\tm)\ni t\mapsto \io \Theta_x^2(\cdot,t)$ ensured
  by (\ref{l1}),
  \bas
	\io \Theta_x^2 
	+ \int_0^t \io \Theta_{xx}^2
	\le \io \Theta_{0x}^2
	+ \int_0^t \io g^2
	\le M^2
	+ c_1
  \eas
  for all $t\in (0,T_0) \cap (0,\tm)$.
\qed
In view of the extensibility criterion in (\ref{ext}), our main result on global existence in the presence of
smooth initial data has thereby been achieved:\abs
\proofc of Theorem \ref{theo8}. \quad
  The statement is an immediate consequence of Lemma \ref{lem5} when combined with Lemma \ref{lem1} and Lemma \ref{lem_loc}.
\qed
\mysection{Further regularity features of $u$ and $v$}\label{sect6}
In this section we prepare our proof of Theorem \ref{theo15} by deriving some further regularity features of the
global solutions just obtained. 
This will be achieved by an analysis of a higher-order quasi-energy functional, the regularity
of which near $t=0$ will result from the following.
\begin{lem}\label{lem55}
  Let $\gamma_0>0$ and $\del(\gamma_0)$ be as in Lemma \ref{lem4}, let $a>0$, and assume (\ref{gf}), (\ref{gg}), (\ref{f})
  and (\ref{init}).
  Then the global solution of (\ref{0}) found in Theorem \ref{theo8} has the additional property that
  \bas
	u \in C^0([0,\infty);W^{2,2}(\Om)).
  \eas
\end{lem}
\proof
  Let $T>0$. Then using that $A:=\gamma(\Theta)$ and $g:=-f'(\Theta)\Theta_x + \gamma'(\Theta) \Theta_x v_x$ satisfiy
  $\{A,A_x,g\} \subset C^0(\bom\times [0,T])$ and $g|_{\pO\times (0,T)}=0$ as well as $A>0$ in $\bom\times [0,T]$,
  in view of (\ref{init}) we can find $c_1>0$ and $c_2>0$, $(v_{0\eta})_{\eta\in (0,1)} \subset C_0^\infty(\Om)$ and
  $(u_{0\eta})_{\eta\in (0,1)} \subset C_0^\infty(\Om)$ as well as $(A_\eta)_{\eta\in (0,1)} \subset C^\infty(\bom\times [0,T])$
  and $(g_\eta)_{\eta\in (0,1)} \subset C_0^\infty(\Om\times [0,T])$ such that
  \be{55.1}
	A_\eta \ge c_1
	\quad \mbox{and} \quad
	|g_\eta| \le c_2
	\quad \mbox{in } \Om\times (0,T)
	\qquad \mbox{for all } \eta\in (0,1),
  \ee
  and that
  \be{55.2}
	\|v_{0\eta} - u_{0t}\|_{W^{1,2}(\Om)}^2
	+ \|u_{0\eta} - u_0\|_{W^{2,2}(\Om)}^2
	\le \eta^\frac{1}{2}
	\qquad \mbox{for all } \eta\in (0,1)
  \ee
  as well as
  \be{55.3}
	\|A_\eta-A\|_{L^\infty(\Om\times (0,T))}^2 
	+ \|g_\eta-g\|_{L^\infty(\Om\times (0,T))}^2 
	\le \eta^\frac{1}{2}
	\qquad \mbox{for all } \eta\in (0,1).
  \ee
  Standard parabolic theory (\cite{amann}, \cite{LSU}) then asserts that for each $\eta\in (0,1)$ there exists $T_\eta\in (0,T]$
  such that the triangular cross-diffusion system
  \be{55.4}
	\lball
	v_{\eta t} = A_\eta(x,t) v_{\eta xx} + au_{\eta xx} + g_\eta(x,t),
	\qquad & x\in\Om, \ t\in (0,T_\eta), \\[1mm]
	u_{\eta t} = \eta u_{\eta xx} + v_\eta,
	\qquad & x\in\Om, \ t\in (0,T_\eta), \\[1mm]
	v_\eta=u_\eta=0,
	\qquad & x\in\pO, \ t\in (0,T_\eta), \\[1mm]
	v_\eta(x,0)=v_{0\eta}(x), \quad u_\eta(x,0)=u_{0\eta}(x),
	\qquad & x\in\Om,
	\ear
  \ee
  admits a classical solution $(v_\eta,u_\eta) \in C^\infty(\bom\times [0,T_\eta))$ which is such that
  \be{55.5}
	\mbox{if $T_\eta<T$, \quad then \quad }
	\limsup_{t\nearrow T_\eta} \Big\{ \|v_\eta(\cdot,t)\|_{W^{1,2}(\Om)} + \|u_\eta(\cdot,t)\|_{W^{1,2}(\Om)} \Big\} =\infty.
  \ee
  On testing the first equation in (\ref{55.4}) against $-v_{\eta xx}$ and using Young's inequality along with (\ref{55.1}),
  we see that these solutions satisfy
  \bas
	\frac{1}{2} \frac{d}{dt} \io v_{\eta x}^2 + \io A_\eta v_{\eta xx}^2
	&=& - a \io u_{\eta xx} (u_{\eta xxt} - \eta u_{\eta xxxx})
	- \io g_\eta v_{\eta xx} \\
	&=& - \frac{a}{2} \frac{d}{dt} \io u_{\eta xx}^2
	- \eta a \io u_{\eta xxx}^2
	- \io g_\eta v_{\eta xx} \\
	&=& - \frac{a}{2} \frac{d}{dt} \io u_{\eta xx}^2
	- \eta a \io u_{\eta xxx}^2
	+ \frac{c_1}{2} \io v_{\eta xx}^2 + \frac{c_2^2 |\Om|}{c_1}
	\qquad \mbox{for all } t\in (0,T_\eta),
  \eas
  because $u_{\eta xx} = u_{\eta t}-v_\eta=0$ on $\pO\times (0,T_\eta)$ by (\ref{55.4}).
  In view of (\ref{55.2}) and the first equation in (\ref{55.1}), this shows that there exist $c_3>0$ and $c_4>0$ fulfilling
  \be{55.55}
	\io v_{\eta x}^2 + \io u_{\eta xx}^2 \le c_3
	\qquad \mbox{for all $t\in (0,T_\eta)$ and } \eta\in (0,1)
  \ee
  as well as
  \be{55.56}
	\int_0^{T_\eta} \io v_{\eta xx}^2 + \eta \int_0^{T_\eta} \io u_{\eta xxx}^2 \le c_4
	\qquad \mbox{for all } \eta\in (0,1).
  \ee
  Due to (\ref{55.5}), this particularly ensures that, in fact, $T_\eta=T$ for all $\eta\in (0,1)$,
  and that with some $(\eta_j)_{j\in\N}\subset (0,1)$, $\wh{v}\in L^2((0,T);W^{2,2}(\Om))$ and
  $\wh{u}\in L^\infty((0,T);W^{2,2}(\Om))$ we have $\eta_j\searrow 0$ as $j\to\infty$ and
  \be{55.6}
	v_\eta \wto \wh{v}
	\quad \mbox{in } L^2((0,T);W^{2,2}(\Om))
	\quad \mbox{and} \quad
	u_\eta \wsto \wh{u}
	\quad \mbox{in } L^\infty((0,T);W^{2,2}(\Om))
  \ee
  as $\eta=\eta_j\searrow 0$.\abs
  On the other hand, a combination of (\ref{55.4}) with (\ref{0v}) shows that
  \bea{55.7}
	\frac{1}{2} \frac{d}{dt} \io (v_\eta-v)^2
	&=& \io (v_\eta-v) \cdot \Big\{ \big( A_\eta v_{\eta xx} + a u_{\eta xx} + g_\eta \big) - \big( Av_{xx} + au_{xx} + g\big)
		\Big\} \nn\\
	&=& \io A(v_\eta-v) (v_\eta-v)_{xx} 
	+ \io (A_\eta-A) (v_\eta-v) v_{\eta xx} \nn\\
	& & + a \io (u_\eta-u)_t (u_\eta-u)_{xx}
	- \eta a \io u_{\eta xx}^2 
	+ \eta a \io u_{\eta xx} u_{xx} \nn\\
	& & + \io (v_\eta-v) (g_\eta-g) \nn\\
	&=& - \io A (v_\eta-v)_x^2
	- \io A_x (v_\eta-v) (v_\eta-v)_x
	+ \io (A_\eta-A) (v_\eta-v)v_{\eta xx} \nn\\
	& & - \frac{a}{2} \frac{d}{dt} \io (u_\eta-u)_x^2
	- \eta a \io u_{\eta xx}^2
	- \eta a \io u_{\eta xxx} u_x \nn\\
	& & + \io (v_\eta-v)(g_\eta-g)
	\qquad \mbox{for all } t\in (0,T),
  \eea
  where
  \bas
	- \io A_x (v_\eta-v) (v_\eta-v)_x
	\le \io A (v_\eta-v)_x^2
	+ \frac{1}{4} \io \frac{A_x^2}{A} (v_\eta-v)^2
  \eas
  and
  \bas
	\io (A_\eta-A) (v_\eta-v)v_{\eta xx} 
	\le \frac{1}{2} \io (v_\eta-v)^2
	+ \frac{1}{2} \|A_\eta-A\|_{L^\infty(\Om)}^2 \io v_{\eta xx}^2
  \eas
  as well as
  \bas
	- \eta a \io u_{\eta xxx} u_x 
	\le \frac{1}{2} \eta^\frac{3}{2} \io u_{\eta xxx}^2 + \frac{\eta^\frac{1}{2} a^2}{2} \io u_x^2
  \eas
  and
  \bas
	\io (v_\eta-v)(g_\eta-g)
	\le \frac{1}{2} \io (v_\eta-v)^2
	+ \frac{|\Om|}{2} \|g_\eta-g\|_{L^\infty(\Om)}^2
  \eas
  for all $t\in (0,T)$ and $\eta\in (0,1)$.
  Since $u_x\in L^\infty((0,T);L^2(\Om))$ by Theorem \ref{theo8}, in line with (\ref{55.3}) and the boundedness of
  $\frac{A_x^2}{A}$ on $\Om\times (0,T)$, we therefore obtain $c_5>0$ such that
  $y_\eta(t):=\io (v_\eta-v)^2(\cdot,t) + a \io (u_\eta-u)_x^2(\cdot,t)$, $t\in [0,T]$, $\eta\in (0,1)$, satisfies
  \bas
	y_\eta'(t) \le c_5 y_\eta(t) + \eta^\frac{1}{2} h_\eta(t)
	\qquad \mbox{for all $t\in (0,T)$ and } \eta\in (0,1),
  \eas
  where (\ref{55.56}) guarantees that for 
  $h_\eta(t):=\io v_{\eta xx}^2(\cdot,t) + \eta \io u_{\eta xxx}^2(\cdot,t) + |\Om|$, $t\in (0,T), \eta\in (0,1)$,
  we have
  \bas
	\int_0^T h_\eta(t) dt \le c_6:=c_4 + |\Om| T
	\qquad \mbox{for all } \eta\in (0,1).
  \eas
  Therefore, 
  \bas
	y_\eta(t)
	\le y_\eta(0) e^{c_5 t} + \eta^\frac{1}{2} \int_0^t e^{c_5(t-s)} h_\eta(s) ds
	\le y_\eta(0) e^{c_5 T} + \eta^\frac{1}{2} c_6 e^{c_5 T}
	\qquad \mbox{for all $t\in (0,T)$ and } \eta\in (0,1),
  \eas
  whence due to (\ref{55.2}) we can find $c_7>0$ such that
  \bas
	y_\eta(t) \le c_7 \eta^\frac{1}{2}
	\qquad \mbox{for all $t\in (0,T)$ and } \eta\in (0,1).
  \eas
  In (\ref{55.6}), we thus must have $\wh{v}=v$ and $\wh{u}=u$ a.e.~in $\Om\times (0,T)$, which especially implies that
  $v=u_t\in L^2((0,T);W^{2,2}(\Om))$ and hence $u\in C^0([0,T];W^{2,2}(\Om))$.
  As $T>0$ was arbitrary, this establishes the claim.
\qed
Again through a variational argument, we can turn the information gained in Lemma \ref{lem5} into improved knowledge
on regularity of the components $u$ and $v$ of the global solutions to (\ref{0v}) 
that will become relevant in our reasoning related to Theorem \ref{theo15}.
\begin{lem}\label{lem6}
  If $\gamma_0>0$, if $\del(\gamma_0)$ is as in Lemma \ref{lem4}, if $a>0$, and if (\ref{gf}), (\ref{gg}) and (\ref{f}) hold,
  then given $M>0$ and $T_0>0$ one can fix $\Lambda_4=\Lambda_4(M,T_0,a,\gamma,f)>0$ 
  with the property that assuming (\ref{init}) and (\ref{iM}), we have
  \be{6.1}
	\io v_x^2(x,t) dx \le \Lambda_4
	\qquad \mbox{for all } t\in (0,T_0)
  \ee
  and
  \be{6.2}
	\io u_{xx}^2(x,t) dx \le \Lambda_4
	\qquad \mbox{for all } t\in (0,T_0)
  \ee
  as well as
  \be{6.3}
	\int_0^T \io v_{xx}^2 \le \Lambda_4
	\qquad \mbox{for all } T\in (0,T_0).
  \ee
\end{lem}
\proof
  We integrate by parts in the first equation from (\ref{0v}) to see that
  \bea{6.33}
	\frac{1}{2} \frac{d}{dt} \io v_x^2
	&=& - \io v_{xx} \cdot \big\{ \gamma(\Theta) v_{xx} + \gamma'(\Theta) \Theta_x v_x + a u_{xx} - f'(\Theta) \Theta_x \big\}
		\nn\\
	&=& - \io \gamma(\Theta) v_{xx}^2 
	- \io \gamma'(\Theta) \Theta_x v_x v_{xx}
	- \frac{a}{2} \frac{d}{dt} \io u_{xx}^2 \nn\\
	& & + \io f'(\Theta) \Theta_x v_{xx}
	\qquad \mbox{for all } t>0,
  \eea
  where by the lower bound in (\ref{gg}),
  \be{6.4}
	- \io \gamma(\Theta) v_x^2 \le - \gamma_0 \io v_{xx}^2
	\qquad \mbox{for all } t>0.
  \ee
  Now taking $c_1>0$ such that in line with the continuity of the embedding $W^{1,2}(\Om) \hra L^\infty(\Om)$ we have
  \bas
	\|\vp\|_{L^\infty(\Om)} \le c_1 \|\vp_x\|_{L^2(\Om)} + c_1 \|\vp\|_{L^1(\Om)}
	\qquad \mbox{for all } \vp \in W^{1,2}(\Om),
  \eas
  from Lemma \ref{lem5} together with Lemma \ref{lem1} we obtain that with $\Lambda_3=\Lambda_3(M,T_0,a,\gamma,f)$ and
  $\Lambda_1=\Lambda_1(M)$ as found there,
  not only
  \be{6.5}
	\io \Theta_x^2 \le \Lambda_3
	\qquad \mbox{for all } t\in (0,T_0) 
  \ee
  but moreover also
  \be{6.6}
	\|\Theta\|_{L^\infty(\Om)}
	\le c_2\equiv c_2(M,T_0,a,\gamma,f):=c_1\Lambda_3^\frac{1}{2} + c_1 \Lambda_1
	\qquad \mbox{for all } t\in (0,T_0);
  \ee
  in particular, using the Gagliardo-Nirenberg inequality to fix $c_3>0$ fulfilling
  \be{6.7}
	\|\vp_x\|_{L^\infty(\Om)}^2 \le c_3 \|\vp_{xx}\|_{L^2(\Om)} \|\vp_x\|_{L^2(\Om)} + c_3 \|\vp_x\|_{L^2(\Om)}^2
	\qquad \mbox{for all } \vp\in W^{2,2}(\Om),
  \ee
  from (\ref{6.6}) we infer that
  \be{6.87}
	|\gamma'(\Theta)|
	\le c_4\equiv c_4(M,T_0,a,\gamma,f) := \|\gamma'\|_{L^\infty((0,c_2))} 
	\qquad \mbox{in } \Om\times (0,T_0),
  \ee
  and that
  \be{6.8}
	|f'(\Theta)| \le c_5\equiv c_5(M,T_0,a,\gamma,f):=\|f'\|_{L^\infty((0,c_2))}
	\qquad \mbox{in } \Om\times (0,T_0).
  \ee
  On the right-hand side of (\ref{6.3}), by utilizing Young's inequality as well as (\ref{6.87}), (\ref{6.5}) and (\ref{6.7})
  we can therefore estimate
  \bea{6.9}
	- \io \gamma'(\Theta) \Theta_x v_x v_{xx}
	&\le& \frac{\gamma_0}{8} \io v_{xx}^2
	+ \frac{2}{\gamma_0} \io \gamma'^2(\Theta) \Theta_x^2 v_x^2 \nn\\
	&\le& \frac{\gamma_0}{8} \io v_{xx}^2
	+ \frac{2 c_4^2}{\gamma_0} \|v_x\|_{L^\infty(\Om)}^2 \io \Theta_x^2 \nn\\
	&\le& \frac{\gamma_0}{8} \io v_{xx}^2
	+ \frac{2 c_4^2 \Lambda_3}{\gamma_0} \|v_x\|_{L^\infty(\Om)}^2 \nn\\
	&\le& \frac{\gamma_0}{8} \io v_{xx}^2
	+ \frac{2 c_3 c_4^2 \Lambda_3}{\gamma_0} \|v_{xx}\|_{L^2(\Om)} \|v_x\|_{L^2(\Om)}
	+ \frac{2 c_3 c_4^2 \Lambda_3}{\gamma_0} \|v_x\|_{L^2(\Om)}^2 \nn\\
	&\le& \frac{\gamma_0}{4} \io v_{xx}^2
	+ c_5 \io v_x^2
	\qquad \mbox{for all } t\in (0,T_0)
  \eea
  with $c_5\equiv c_5(M,T_0,a,\gamma,f):=\frac{8 c_3^2 c_4^4 \Lambda_3^2}{\gamma_0^3} + \frac{2 c_3 c_4^2 \Lambda_3}{\gamma_0}$.
  Apart from that, Young's inequality, (\ref{6.8}) and (\ref{6.5}) ensure that
  \bas
	\io f'(\Theta) \Theta_x v_{xx}
	&\le& \frac{\gamma_0}{4} \io v_{xx}^2
	+ \frac{1}{\gamma_0} \io f'^2(\Theta) \Theta_x^2 \\
	&\le& \frac{\gamma_0}{4} \io v_{xx}^2
	+ \frac{c_5^2}{\gamma_0} \io \Theta_x^2 \\
	&\le& \frac{\gamma_0}{4} \io v_{xx}^2
	+ c_6
	\qquad \mbox{for all } t\in (0,T_0)
  \eas
  with $c_6\equiv c_6(M,T_0,a,\gamma,f):=\frac{c_5^2 \Lambda_3}{\gamma_0}$, so that from (\ref{6.3}), (\ref{6.4}) and (\ref{6.9})
  we altogether conclude that
  \bas
	\frac{d}{dt} \io v_x^2 + a \frac{d}{dt} \io u_{xx}^2
	+ \gamma_0 \io v_{xx}^2
	\le 2c_5 \io v_x^2 + 2c_6
	\qquad \mbox{for all } t\in (0,T_0).
  \eas
  Writing $c_7\equiv c_7(M,T_0,a,\gamma,f):=\max\{2c_5 \, , \, 2c_6\}$, we thus obtain that
  \be{6.10}
	\frac{d}{dt} \bigg\{ 1 + \io v_x^2 + a \io u_{xx}^2 \bigg\}
	+ \gamma_0 \io v_{xx}^2
	\le c_7 \cdot \bigg\{ 1 + \io v_x^2 + a \io u_{xx}^2 \bigg\}
	\qquad \mbox{for all } t\in (0,T_0),
  \ee
  which through a comparison argument firstly entails that, by the Cauchy-Schwarz inequality and (\ref{iM}),
  \bas
	1+ \io v_x^2 + a \io u_{xx}^2 
	&\le& \bigg\{ 1 + \io (u_{0t})_x^2 + a \io u_{0xx}^2 \bigg\} \cdot e^{c_7 T_0} \\
	&\le&
	\Bigg\{ 1 + |\Om|^\frac{1}{2} \bigg\{ \io (u_{0t})_x^4 \bigg\}^\frac{1}{2}  + a \io u_{0xx}^2 \Bigg\} \cdot e^{c_7 T_0} \\
	&\le& c_8\equiv c_8(M,T_0,a,\gamma,f):= 
	( 1 + |\Om|^\frac{1}{2} M^2 + a M^2) \cdot e^{c_7 T_0}
  \eas
  for all $t\in [0,T_0)$,
  because both $[0,\infty) \ni t\mapsto \io v_x^2(\cdot,t)$ and $[0,\infty)\ni t\mapsto \io u_{xx}^2(\cdot,t)$ 
  are continuous according to Theorem \ref{theo8} and Lemma \ref{lem55}.
  Thereafter, upon an integration the inequality in (\ref{6.10}) secondly implies that
  \bas
	\gamma_0 \int_0^T \io v_{xx}^2 \le c_8 + c_7 c_8 T_0
	\qquad \mbox{for all } t\in (0,T_0),
  \eas
  so that with an obvious choice of $\Lambda_4(M,T_0,a,\gamma,f)$, we can achieve (\ref{6.1})-(\ref{6.3}).
\qed
\mysection{Coping with rough initial data. Proof of Theorem \ref{theo15}}\label{sect7}
In order to appropriately cope with the setting of less regular initial data addressed in Theorem \ref{theo15}, 
given $(u_0,u_{0t},\Theta_0)$ fulfilling (\ref{Init}) we fix approximations thereof which are such that
\be{ie}
	\lbal
	(v_{0\eps})_{\eps\in (0,1)} \subset W^{2,2}(\Om) \cap W_0^{1,2}(\Om), \qquad \mbox{that}\\[1mm]
	(u_{0\eps})_{\eps\in (0,1)} \subset W^{3,2}(\Om) \cap W_0^{1,2}(\Om),
	\qquad \mbox{and that} \\[1mm]
	(\Theta_{0\eps})_{\eps\in (0,1)} \subset W^{2,2}(\Om) 
	\mbox{ satisfies $\Theta_{0\eps}\ge 0$ and $\Theta_{0\eps x}|_{\pO}=0$ for all $\eps\in (0,1)$,}
	\ear
\ee
and for which we have
\be{conv}
	v_{0\eps}\to u_{0t}
	\ \mbox{in } W^{1,4}(\Om),
	\quad
	u_{0\eps} \to u_0	
	\ \mbox{in } W^{2,2}(\Om),
	\quad \mbox{and} \quad
	\Theta_{0\eps} \to \Theta_0
	\ \mbox{in } W^{1,2}(\Om)
	\qquad \mbox{as } \eps\searrow 0.
\ee
For $\eps\in (0,1)$, we consider the corresponding version of (\ref{0v}) given by	
\be{0eps}
	\lball
	\vepst = \big( \gamma(\Teps) \vepsx \big)_x + a\uepsxx - \big( f(\Teps)\big)_x,
	\qquad & x\in\Om, \ t>0, \\[1mm]
	\uepst = \veps,
	\qquad & x\in\Om, \ t>0, \\[1mm]
	\Tepst = \Tepsxx + \gamma(\Teps) \vepsx^2 - f(\Teps) \vepsx,
	\qquad & x\in\Om, \ t>0, \\[1mm]
	\veps=0, \quad \ueps=0, \quad \Tepsx=0,
	\qquad & x\in\pO, \ t>0, \\[1mm]
	\veps(x,0)=v_{0\eps}(x), \quad \ueps(x,0)=u_{0\eps}(x), \quad \Teps(x,0)=\Theta_{0\eps}(x),
	\qquad & x\in\Om,
	\ear
\ee
and note that our previously obtained results ensure the following.
\begin{lem}\label{lem11}
  Let $\gamma_0>0$ and $\del=\del(\gamma_0)>0$ be as in Theorem \ref{theo8}, let $a>0$, 
  and assume (\ref{gf}), (\ref{gg}), (\ref{f}) and (\ref{ie}). Then 
  for each $\eps\in (0,1)$, there exist unique functions
  \bas
	\lbal
	\veps \in \Big( \bigcup_{\beta\in (0,1)} C^{1+\beta,\frac{1+\beta}{2}}(\bom\times [0,\infty))\Big) 
		\cap C^{2,1}(\bom\times (0,\infty)), \\[1mm]
	\ueps\in \Big( \bigcup_{\beta\in (0,1)} C^{1+\beta,\frac{1+\beta}{2}}(\bom\times [0,\infty))\Big) 
		\cap C^{2,1}(\bom\times (0,\infty)) \cap C^0([0,\infty);W^{2,2}(\Om))
		\qquad \mbox{and} \\[1mm]
	\Teps\in \Big( \bigcup_{\beta\in (0,1)} C^{1+\beta,\frac{1+\beta}{2}}(\bom\times [0,\infty))\Big) 
		\cap C^{2,1}(\bom\times (0,\infty))
	\ear
  \eas
  such that $\Teps\ge 0$ in $\Om\times [0,\infty)$, and that (\ref{0eps}) is solved in the classical sense.
\end{lem}
\proof
  This actually is a by-product of Theorem \ref{theo8} when combined with Lemma \ref{lem55}.
\qed
As (\ref{conv}) ensures validity of (\ref{iM}) with some $M>0$ independent of suitably small $\eps$,
we can now take advantage of the fact that in their most substantial part our previously gained estimates
depend on the initial data only through this bound $M$.
The following lemma thus actually collects knowledge gathered above only.
\begin{lem}\label{lem12}
  If $\gamma_0>0$ and $\del=\del(\gamma_0)>0$ is as in Theorem \ref{theo8}, let $a>0$, 
  and if (\ref{gf}), (\ref{gg}) and (\ref{f}) and (\ref{Init}) hold,
  then for any $T_0>0$ one can find 
  $C=C(T_0,a,\gamma,f,u_0,u_{0t},\Theta_0)>0$ with the property that if (\ref{ie}) and (\ref{conv}) hold,
  then the solutions of (\ref{0eps}) satisfy
  \be{12.1}
	\io \vepsx^2(\cdot,t) \le C
	\qquad \mbox{for all $t\in (0,T)$ and } \eps\in (0,\eps_0)
  \ee
  and
  \be{12.2}
	\io \uepsxx^2(\cdot,t) \le C
	\qquad \mbox{for all $t\in (0,T)$ and } \eps\in (0,\eps_0)
  \ee
  and
  \be{12.3}
	\|\Teps(\cdot,t)\|_{L^\infty(\Om)} + \io \Tepsx^2(\cdot,t) \le C
	\qquad \mbox{for all $t\in (0,T)$ and } \eps\in (0,\eps_0),
  \ee
  and that
  \be{12.4}
	\int_0^{T_0} \io \vepsxx^2 \le C
	\qquad \mbox{for all } \eps\in (0,\eps_0)
  \ee
  and
  \be{12.5}
	\int_0^{T_0} \io \vepsx^4 \le C
	\qquad \mbox{for all } \eps\in (0,\eps_0)
  \ee
  as well as
  \be{12.6}
	\int_0^{T_0} \io \Tepsxx^2 \le C
	\qquad \mbox{for all } \eps\in (0,\eps_0).
  \ee
\end{lem}
\proof
  Writing
  \bas
	M\equiv M(u_0,u_{0t},\Theta_0)
	&:=& 1+
	\|u_0\|_{L^\infty(\Om)} 
	+ \|u_{0x}\|_{L^\infty(\Om)}
	+ \|u_{0xx}\|_{L^2(\Om)} \\
	& & + \|u_{0t}\|_{L^\infty(\Om)}
	+ \big\| \big(u_{0t}\big)_x\big\|_{L^4(\Om)}
	+ \|\Theta_0\|_{L^\infty(\Om)}
	+ \|\Theta_{0x}\|_{L^2(\Om)}
  \eas
  in view of (\ref{conv}) we can find $\eps_0=\eps_0(u_0,u_{0t},\Theta_0)\in (0,1)$ such that
  \bas
	& & \hs{-20mm}
	\|u_{0\eps}\|_{L^\infty(\Om)} 
	+ \|u_{0\eps x}\|_{L^\infty(\Om)}
	+ \|u_{0\eps xx}\|_{L^2(\Om)} \\
	& & + \|v_{0\eps}\|_{L^\infty(\Om)}
	+ \|v_{0\eps x}\|_{L^4(\Om)}
	+ \|\Theta_{0\eps}\|_{L^\infty(\Om)}
	+ \|\Theta_{0\eps x}\|_{L^2(\Om)}
	\le M
	\qquad \mbox{for all } \eps\in (0,\eps_0).
  \eas
  We may therefore apply Lemma \ref{lem6} to directly obtain (\ref{12.1}), (\ref{12.2}) and (\ref{12.4}),
  while (\ref{12.5}) has been derived in Lemma \ref{lem4}.
  In view of the continuity of the embedding $W^{1,2}(\Om) \hra L^\infty(\Om)$, the estimate in
  (\ref{12.3}) and (\ref{12.6}) results from Lemma \ref{lem5} when combined with Lemma \ref{lem1}.
\qed
These properties imply $L^2$ bounds for the time derivatives $\vepst$ and $\Tepst$:
\begin{lem}\label{lem13}
  Let $\gamma_0>0$ and $\del=\del(\gamma_0)>0$ be taken from Theorem \ref{theo8}, let $a>0$, 
  and suppose that (\ref{gf}), (\ref{gg}) and (\ref{f}) and (\ref{Init}) hold.
  Then for each $T_0>0$ there exists  
  $C=C(T_0,a,\gamma,f,u_0,u_{0t},\Theta_0)>0$ such that 
  whenever (\ref{ie}) and (\ref{conv}) are fulfilled,
  for the solutions of (\ref{0eps}) we have
  \be{13.1}
	\int_0^{T_0} \io \vepst^2 \le C
	\qquad \mbox{for all } \eps\in (0,\eps_0)
  \ee
  and
  \be{13.2}
	\int_0^{T_0} \io \Tepst^2 \le C
	\qquad \mbox{for all } \eps\in (0,\eps_0),
  \ee
  where $\eps_0=\eps_0(u_0,u_{0t},\Theta_0)$ is as in Lemma \ref{lem12}.
\end{lem}
\proof
  Utilizing (\ref{0eps}) and Young's inequality, we can estimate
  \bas
	\int_0^{T_0} \io \vepst^2
	&=& \int_0^{T_0} \io \big\{ \gamma(\Teps) \vepsxx + \gamma'(\Teps)\Tepsx \vepsx + a\uepsxx - f'(\Teps) \Tepsx \big\} \\
	&\le& 4 \int_0^{T_0} \io \gamma^2(\Teps) \vepsxx^2
	+ 4 \int_0^{T_0} \io \gamma'^2(\Teps) \Tepsx^2 \vepsx^2 \\
	& & + 4 a^2 \int_0^{T_0} \io \uepsxx^2
	+ 4 \int_0^{T_0} \io \feps'^2(\Teps) \Tepsx^2
	\qquad \mbox{for all } \eps\in (0,1),
  \eas
  whence using the the continuity of $\gamma, \gamma'$ and $f'$ along with the
  $L^\infty$ bound contained in (\ref{12.3}) we can find $c_1=c_1(T_0,a,\gamma,f,u_0,u_{0t},\Theta_0)>0$ such that
  \bas
	\int_0^{T_0} \io \vepst^2
	\le c_1 \int_0^{T_0} \io \vepsxx^2
	+ c_1 \int_0^{T_0} \io \Tepsx^2 \vepsx^2
	+ c_1 \int_0^{T_0} \io \uepsxx^2
	+ c_1 \int_0^{T_0} \io \Tepsx^2
  \eas
  for all $\eps\in (0,\eps_0)$.
  Here, using that $\Tepsx=0$ on $\pO\times (0,\infty)$ for all $\eps\in (0,1)$,
  and that $\|\vp\|_{L^\infty(\Om)} \le |\Om|^\frac{1}{2} \|\vp_x\|_{L^2(\Om)}$ for all $\vp\in W_0^{1,2}(\Om)$, we see that
  \bas
	\int_0^{T_0} \io \Tepsx^2 \vepsx^2
	\le \int_0^{T_0} \|\Tepsx\|_{L^\infty(\Om)}^2 \|\vepsx\|_{L^2(\Om)}^2
	\le |\Om| \cdot \Big\{ \sup_{t\in (0,T_0)} \|\vepsx(\cdot,t)\|_{L^2(\Om)}^2 \Big\} \cdot \int_0^{T_0} \io \Tepsxx^2
  \eas
  for all $\eps\in (0,1)$.
  The estimate in (\ref{13.1}) therefore follows upon combining (\ref{12.4}) and (\ref{12.1}) with (\ref{12.6}), (\ref{12.2}) and 
  (\ref{12.3}).\abs
  Likewise, from the third equation in (\ref{0eps}), (\ref{12.3}) and the local boundedness of $\gamma$ and $f$ in 
  $[0,\infty)$ we obtain that with some $c_2=c_2(T_0,a,\gamma,f,u_0,u_{0t},\Theta_0)>0$,
  \bas
	\int_0^{T_0} \io \Tepst^2
	&\le& 3 \int_0^{T_0} \io \Tepsxx^2
	+ 3 \int_0^{T_0} \io \gamma^2(\Teps) \vepsx^4
	+ 3 \io f^2(\Teps) \vepsx^2 \\
	&\le& c_2 \int_0^{T_0} \io \Tepsxx^2
	+ c_2 \int_0^{T_0} \io \vepsx^4
	+ c_2 \int_0^{T_0} \io \vepsx^2
  \eas
  for all $\eps\in (0,\eps_0)$. 
  From (\ref{12.6}) and (\ref{12.5}) we thus infer that 
  upon enlarging $C$ if necessary we can also achieve (\ref{13.2}).
\qed
In a straightforward manner, Lemma \ref{lem12} and Lemma \ref{lem13} facilitate 
the construction of a limit which solves (\ref{0}) in the sense specified in Theorem \ref{theo15}:
\begin{lem}\label{lem14}
  Let $\gamma_0>0$, let $\del=\del(\gamma_0)>0$ be as in Theorem \ref{theo8}, let $a>0$, 
  and assume (\ref{gf}), (\ref{gg}), (\ref{f}), (\ref{Init}), (\ref{ie}) and (\ref{conv}).
  Then there exist $(\eps_j)_{j\in\N} \subset (0,1)$ as well as
  \be{14.1}
	\lbal
	v\in C^0(\bom\times [0,\infty)) \cap L^\infty_{loc}((0,\infty);W_0^{1,2}(\Om)) 
		\cap L^2_{loc}([0,\infty);W^{2,2}(\Om)), \\[1mm]
	u\in C^0(\bom\times [0,\infty)) \cap L^\infty_{loc}((0,\infty);W^{2,2}(\Om) \cap W_0^{1,2}(\Om))
	\qquad \mbox{and} \qquad \\[1mm]
	\Theta\in C^0(\bom\times [0,\infty)) \cap L^\infty_{loc}([0,\infty);W^{1,2}(\Om)) \cap L^2_{loc}([0,\infty);W_N^{2,2}(\Om))
	\ear
  \ee
  such that 
  \be{14.01}
	v(\cdot,0)=u_{0t},
	\quad
	u(\cdot,0)=u_0
	\quad \mbox{and} \quad
	\Theta(\cdot,0)=\Theta_0
	\qquad \mbox{in } \Om,
  \ee
  that $\Theta\ge 0$ in $\Om\times (0,\infty)$, that $\eps_j\searrow 0$ as $j\to\infty$,
  and that with $((\veps,\ueps,\Teps))_{\eps\in (0,1)}$ taken from Lemma \ref{lem11} we have
  \begin{eqnarray}
	& & \veps\to v
	\qquad \mbox{in } C^0_{loc}(\bom\times [0,\infty)),
	\label{14.2} \\
	& & \vepsx \to v_x
	\qquad \mbox{a.e.~in $\Om\times (0,\infty)$ and in } L^2_{loc}(\bom\times [0,\infty)),
	\label{14.3} \\
	& & \vepsxx \wto v_{xx}
	\qquad \mbox{ in } L^2_{loc}(\bom\times [0,\infty)),
	\label{14.31} \\
	& & \vepst \wto v_t
	\qquad \mbox{ in } L^2_{loc}(\bom\times [0,\infty)),
	\label{14.32} \\
	& & \ueps \to u
	\qquad \mbox{in } C^0_{loc}(\bom\times [0,\infty)),
	\label{14.4} \\
	& & \uepsxx \wto u_{xx}
	\qquad \mbox{in } L^2_{loc}(\bom\times [0,\infty)),
	\label{14.5} \\
	& & \uepst \wto u_t
	\qquad \mbox{in } L^2_{loc}(\bom\times [0,\infty)),
	\label{14.51} \\
	& & \Teps \to \Theta
	\qquad \mbox{in } C^0_{loc}(\bom\times [0,\infty)),
	\label{14.6} \\
	& & \Tepsx \to \Theta_x
	\qquad \mbox{a.e.~in $\Om\times (0,\infty)$ and in } L^2_{loc}(\bom\times [0,\infty)),
	\label{14.7} \\
	& & \Tepsxx \wto \Theta_{xx}
	\qquad \mbox{ in } L^2_{loc}(\bom\times [0,\infty))
	\qquad \mbox{and}
	\label{14.8} \\
	& & \Tepst \wto \Theta_t
	\qquad \mbox{in } L^2_{loc}(\bom\times [0,\infty))
	\label{14.9} 
  \end{eqnarray}
  as $\eps=\eps_j\searrow 0$.
  Moreover, (\ref{wt}) as well as the identities $u_t=v$ and
  \be{wv}
	v_t = \gamma(\Theta) v_{xx} + \gamma'(\Theta) \Theta_x v_x + a u_{xx} - f'(\Theta)\Theta_x
  \ee
  hold a.e.~in $\Om\times (0,\infty)$.
\end{lem}
\proof
  For each $T_0>0$, from Lemma \ref{lem12} we know that 
  \bas
	(\veps)_{\eps\in (0,1)} 
	\ \mbox{ is bounded in $C^0(([0,T_0];W^{1,2}(\Om))$ and in } L^2((0,T_0);W^{2,2}(\Om)),
  \eas
  that
  \bas
	(\ueps)_{\eps\in (0,1)} 
	\ \mbox{ is bounded in } L^\infty((0,T_0);W^{2,2}(\Om)),
  \eas
  and that
  \bas
	(\Teps)_{\eps\in (0,1)} 
	\ \mbox{ is bounded in $C^0([0,T_0];W^{1,2}(\Om))$ and in } L^2((0,T_0);W^{2,2}(\Om)),
  \eas
  while Lemma \ref{lem13} together with the fact that $\uepst=\veps$ warrants that
  \bas
	(\vepst)_{\eps\in (0,1)},
	\ 
	(\uepst)_{\eps\in (0,1)}
	\ \mbox{and} \ 
	(\Tepst)_{\eps\in (0,1)}
	\ \mbox{ are bounded in } L^2(\Om\times (0,T_0)).
  \eas
  Thanks to the compactness of the embeddings from $W^{1,2}(\Om)$ into $C^0(\bom)$ and from $W^{2,2}(\Om)$ into $W^{1,2}(\Om)$,
  we may thus employ an Aubin-Lions lemma to find $(\eps_j)_{j\in\N}\subset (0,1)$ such that $\eps_j\searrow 0$
  as $j\to\infty$, and that (\ref{14.2})-(\ref{14.9}) hold with some functions $v,u$ and $\Theta$
  fulfilling (\ref{14.1}). While (\ref{14.01}) then immediately results from (\ref{0eps}), (\ref{14.2}), (\ref{14.4}) and
  (\ref{14.6}),
  the solution properties (\ref{wv}) and (\ref{wt}) can readily be derived from (\ref{14.3})-(\ref{14.32}) and 
  (\ref{14.5})-(\ref{14.9}) by taking limits in the identities
  \bas
	\int_0^\infty \io \uepst \vp = \int_0^\infty \io \veps \vp
  \eas
  and
  \bas
	\int_0^\infty \io \vepst \vp = \int_0^{T_0} \io \gamma(\Teps)\vepsxx \vp
	+ \int_0^\infty \io \gamma'(\Teps)\Tepsx \vepsx \vp
	+ a \int_0^\infty \io \uepsxx \vp
	- \int_0^\infty \io f'(\Teps) \Tepsx \vp
  \eas
  and
  \bas
	\int_0^\infty \io \Tepst \vp = \int_0^{T_0} \io \Tepsxx \vp
	+ \int_0^\infty \io \gamma(\Teps) \vepsx^2 \vp
	- \int_0^\infty \io f(\Teps) \vepsx \vp,
  \eas
  as known to be valid for each $\eps\in (0,1)$ and any $\vp\in C_0^\infty(\Om\times (0,\infty))$ by (\ref{0eps}).
\qed
Our main result on global strong solvability in (\ref{0}) has thereby been accomplished:\abs
\proofc of Theorem \ref{theo15}.\quad
  All statements immediately result from Lemma \ref{lem14}.
\qed

\bigskip

{\bf Acknowlegements.} \quad
The author acknowledges support of the Deutsche Forschungsgemeinschaft (Project No. 444955436).\abs
{\bf Conflict of interest statement.} \quad
The author declares that he has no conflict of interest, 
and that he has no relevant financial or non-financial interests to disclose.\abs
{\bf Data availability statement.} \quad
Data sharing is not applicable to this article as no datasets were
generated or analyzed during the current study.\abs

\end{document}